\documentclass[11pt]{amsart}
\usepackage{amsmath,amstext,amssymb,amsopn,amsthm,mathrsfs}
\newtheorem{theo}{Theorem}[section]
\newtheorem{prop}{Proposition}[section]

\numberwithin{equation}{section}
\newcommand{\dem}{\medskip \par \noindent \mbox{\bf Proof }}
\def\ep{\rule{1.5mm}{3mm}}
\begin{document}

\title[Transference result Jacobi Riesz transform]{A transference result of the $L^p$ continuity of the Jacobi Riesz transform to the Gaussian and Laguerre Riesz transforms.}

\author{Eduard Navas}
\address{Departmento de Matem\'aticas, Universidad Experimental Francisco de Miranda, Punto Fijo, Venezuela.}
\email{[Eduard Navas]eduard.navas@gmail.com}
\author{Wilfredo O. Urbina}
\address{
 Department of Mathematical and Actuarial Sciences, Roosevelt University  Chicago, Il, 60605, USA.}
\email{[Wilfredo Urbina]wurbinaromero@roosevelt.edu}

\begin{abstract}
In this paper using the well known asymptotic relations between 
Jacobi polynomials and Hermite and Laguerre polynomials.
we develop a transference method  to obtain the $L^p$-continuity of the Gaussian-Riesz transform and the $L^p$-continuity of the Laguerre-Riesz transform
from the $L^p$-continuity of the Jacobi-Riesz transform, in dimension one.
 \end{abstract}

\maketitle

\begin{center}
\large{\em Dedicated to professor  Calixto P. Calder\'on
with deep admiration.}
\end{center}

\section{Preliminaries}
 In the theory of classical orthogonal polynomials, it  is well know the asymp-totic relations between  Jacobi
polynomials and  Hermite and Laguerre polynomials.
Base in those asymptotic relations  we develop a   transference method to obtain  $L^p$-continuity for the Gaussian-Riesz transform and the $L^p$-continuity of the Laguerre-Riesz transform
from the $L^p$-continuity of the Jacobi-Riesz transform, in dimension one.

 For all the classical polynomials we are going to use the normalizations given  in G. Seg\"o's book \cite{sz}.

\begin{itemize}

\item {\em Jacobi polynomials:}  
For  $\alpha,\beta > -1,$ the Jacobi polynomials $\{P^{(\alpha,\beta)}_n\}_{n \in \mathbb N }$ are defined
as the orthogonal  polynomials associated with the Jacobi measure $\mu_{\alpha,\beta}$
(or beta measure) in $(-1,1)$, defined as
 
 \begin{eqnarray}
\mu_{\alpha,\beta}(dx) &=& \omega_{\alpha,\beta}(x) dx =
\chi_{(-1,1)}(x)\frac{(1-x)^{\alpha}(1+x)^{\beta}}{2^{\alpha+\beta+1} B(\alpha+1,\beta+1)}dx\\
&=& \nonumber \eta_{\alpha,\beta}
\chi_{(-1,1)}(x)(1-x)^{\alpha}(1+x)^{\beta} dx,
\end{eqnarray}
where $ \eta_{\alpha,\beta} = \frac{1}{2^{\alpha+\beta+1}
B(\alpha+1,\beta+1)} =
\frac{\Gamma(\alpha+\beta+2)}{2^{\alpha+\beta+1}
\Gamma(\alpha+1)\Gamma(\beta+1)}$.

The function  $\omega_{\alpha,\beta}$ is called the (normalized)
{\em Jacobi weight}.\\

The Jacobi polynomials can be obtained from the canonical
basis of the polynomials $\{ 1, x, x^2,\cdots,x^n,\cdots \} $
using the Gram-Schmidt orthogonalization process with respect to
the inner product in $L^2(\mu_{\alpha,\beta})$. Thus we have the {\em
orthogonality property} of Jacobi polynomials  with respect
to  $\mu_{\alpha,\beta} $, 
\begin{equation} \label{proJacOrt}
\int^{\infty}_{-\infty} P^{(\alpha,\beta)}_n(y)
P^{(\alpha,\beta)}_m(y) \, \mu_{\alpha,\beta}(dy)
=\eta_{\alpha,\beta}  h_{n}^{\left( \alpha ,\beta
\right)}\delta_{n,m} = \hat{h_{n}}^{\left( \alpha ,\beta
\right)}\delta_{n,m},
\end{equation}
$n,m =0,1,2, \cdots$, where
\begin{equation}
h_{n}^{\left( \alpha ,\beta \right)}=\frac{ 2^{\alpha +\beta
+1}}{(2n+\alpha +\beta +1)}\frac{\Gamma \left( n+\alpha +1\right)
\Gamma \left( n+\beta +1\right) }{\Gamma \left( n+1\right) \Gamma
\left( n+\alpha +\beta +1\right) },
\end{equation}
and
\begin{eqnarray*}\label{l2norm}
 \hat{h_{n}}^{\left( \alpha ,\beta \right)}&=&\frac{1}{(2n+\alpha
+\beta +1)}\frac{\Gamma(\alpha+\beta+2)\Gamma \left( n+\alpha
+1\right) \Gamma \left( n+\beta +1\right)
}{ \Gamma(\alpha+1)\Gamma(\beta+1)\Gamma \left( n+1\right) \Gamma \left( n+\alpha +\beta +1\right) }\\
&=&  \| P_n^{\left( \alpha ,\beta\right)}\|_{2,(\alpha,\beta)}^{2}.
\end{eqnarray*}
Moreover,

\begin{eqnarray}\label{normPn-1}
 \nonumber&&\left\|P^{(\alpha+1,\beta+1)}_{n-1}\right\|^{2}_{2,(\alpha+1,\beta+1)}\\
\nonumber&=&\frac{1}{(2n+\alpha+\beta+1)}\frac{\Gamma(\alpha+\beta+4)\Gamma
\left( n+\alpha+1\right) \Gamma \left( n+\beta+1\right)
}{ \Gamma(\alpha+2)\Gamma(\beta+2)\Gamma \left( n\right) \Gamma \left( n+\alpha+\beta+2\right) }\\
\nonumber&=&\frac{(\alpha+\beta+3)(\alpha+\beta+2)n }{(\alpha+1)(\beta+1)
 \left( n+\alpha+\beta+1\right)}\\
&&\quad \quad \times\frac{\Gamma(\alpha+\beta+2)\Gamma\left(
n+\alpha+1\right) \Gamma
\left(n+\beta+1\right)}{(2n+\alpha+\beta+1)\Gamma(\alpha+1)\Gamma(\beta+1)\Gamma
\left( n+1\right) \Gamma \left(n+\alpha+\beta+1\right)}\\
\nonumber&=&\frac{(\alpha+\beta+3)(\alpha+\beta+2)n }{(\alpha+1)(\beta+1)
 \left( n+\alpha+\beta+1\right)}\left\|P^{(\alpha,\beta)}_{n}\right\|^{2}_{2,(\alpha,\beta)}.\\\nonumber 
\end{eqnarray}

Now using  the generalized Rodrigues' formula,  Szeg{\"o} \cite{sz},(4.10.1), for $m=1,$ we get

\begin{equation}\label{proPn-1}
\frac{d}{dx}\left\{(1-x)^{\alpha+1}(1+x)^{\beta+1} P_{n-1}^{\left(
\alpha+1,\beta+1\right) }\left( x\right) \right\}
=-2n(1-x)^{\alpha}(1+x)^{\beta}P_{n}^{\left( \alpha,\beta\right)
}\left(x\right).
\end{equation}

On the other hand, the Jacobi polynomial $P_n^{\left( \alpha ,\beta\right)}$ is a
polynomial solution of the {\em Jacobi differential equation},
with parameters $\alpha ,\beta, n$,
\begin{equation}\label{ecuadif}
\left( 1-x^{2}\right) y^{\prime \prime }+\left[ \beta -\alpha
-\left( \alpha +\beta +2\right) x\right] y^{\prime }+n\left(
n+\alpha +\beta +1\right) y=0,
\end{equation}
i. e.  $P_n^{\left( \alpha ,\beta\right)}$ is an eigenfunction of
the (one-dimensional) second order diffusion operator
\begin{equation}\label{jacobiop}
{\mathcal L}^{\alpha,\beta}= (1-x^2)  \frac{d^2}{dx^2} +(\beta-\alpha-(\alpha+\beta+2)x)\frac{d}{dx},
\end{equation}
associated with the eigenvalue $- \lambda^{\alpha+\beta}_n= n(n+\alpha+\beta+1)$. ${\mathcal L}^{\alpha,\beta}$ is  called  the {\em Jacobi differential operator }.\\

The operator semigroup associated to the Jacobi polynomials is
defined for positive or bounded measurable Borel functions of
 $(-1,1)$, as
\begin{equation}\label{Jac1}
T^{\alpha,\beta}_t f(x)=    \int_{-1}^{1}
p^{\alpha,\beta}(t,x,y)f(y) \mu_{\alpha,\beta}(dy),
\end{equation}
where
$$p^{\alpha,\beta}(t,x,y) =   \sum_k \frac{e^{- k(k+\alpha+\beta+1) t}}{ \hat{h_{k}}^{\left( \alpha ,\beta \right)}} P^{(\alpha,\beta)}_k(x) P^{(\alpha,\beta)}_k(y).$$
The  explicit representation of $p^{\alpha,\beta}(t,x,y)$ is
very complicated since  the eigenvalues $\lambda_n$  are not
linearly distributed and was obtained by G. Gasper \cite{gasp}. It is an analog of Bailey's $F_4$ representation for what is usually called the Jacobi-Poisson integral, see \cite{bai}. Considering the case From that form, taking $x=-y=1$, it can be proved that $p^{\alpha,\beta}(t, x, y) $ is a positive kernel for any $x\in (-1,1)$.\\

$\{ T^{\alpha,\beta}_t \}$ is called  the  {\em  Jacobi semigroup}
and can be proved that is a Markov semigroup, for details see
\cite{ur2}.

The Jacobi-Poisson semigroup $\{ P^{\alpha,\beta}_t \}$ can be
defined, using Bochner's {subordination} formula,
\begin{eqnarray*}
e^{-\lambda^{1/2}t} = \frac{1}{\sqrt{\pi}} \int_0^\infty \frac{e^{-u}
}{\sqrt{u}} e^{-{\frac{\lambda t^2}{4u}}}du.
\end{eqnarray*}
as the subordinated semigroup of the Jacobi semigroup,
\begin{eqnarray*}
P_t^{\alpha,\beta} f(x) = \frac{1}{\sqrt \pi} \int_0^{\infty}
\frac{e^{-u}}{\sqrt u} T^{\alpha,\beta}_{t^2/4u}f(x) du.
\end{eqnarray*}

For a function $f\in
L^{2}\left(\left[-1,1\right],\mu_{(\alpha,\beta)}\right)$ let us
consider its Fourier- Jacobi expansion
\begin{equation}\label{jacobides}
f= \sum _{k=0}^\infty\frac{\langle f,P_{k}^{(\alpha,\beta)}
\rangle}{ \hat{h_{k}}^{\left( \alpha ,\beta \right)}}
P_{k}^{(\alpha,\beta)},
\end{equation}
where
$$\langle f,P_{k}^{(\alpha,\beta)} \rangle= \int_{-1}^{1}f(y) P_{k}^{(\alpha,\beta)}(y) \mu_{\alpha,\beta}(dy).$$
 Then, the action of $T_t$ and $P_t$ can be expressed as
\begin{equation*}
T_t^{\alpha,\beta} f =  \sum _{k=0}^\infty \frac{\langle f,P_{k}^{(\alpha,\beta)}
\rangle}{ \hat{h_{k}}^{\left( \alpha ,\beta \right)}}
e^{-\lambda_k t} P_{k}^{(\alpha,\beta)},
\end{equation*}
and
\begin{equation*}
P_t^{\alpha,\beta} f = \sum _{k=0}^\infty \frac{\langle f,P_{k}^{(\alpha,\beta)}
\rangle}{ \hat{h_{k}}^{\left( \alpha ,\beta \right)}}
e^{-\sqrt{\lambda_k} t} P_{k}^{(\alpha,\beta)}.
\end{equation*}
Following  the classical case, see \cite{noStemp},  the {\em Jacobi-Riesz transform} can be
define formally as
 \begin{equation}\label{rieszJac}
R^{\alpha, \beta} =  \sqrt{1-x^2} \frac{d}{dx} ({\mathcal L}^{\alpha,\beta})^{-1/2},
\end{equation}
where $({\mathcal L}^{\alpha,\beta})^{- \nu/2}$ is the
Jacobi-Riesz potential of order $\nu/2$.  $({\mathcal L}^{\alpha,\beta})^{- \nu/2}$ can be represented
as
$$ ({\mathcal L}^{\alpha,\beta})^{- \nu /2} f =\frac{1}{\Gamma(\nu)} \int_0^\infty t^{\nu-1}
P_t^{(\alpha,\beta)}  f dt,$$ and then for
$f\in L^{2}\left(\left[-1,1\right],\mu_{(\alpha,\beta)}\right)$ ,
$({\mathcal L}^{\alpha,\beta})^{- \nu/2}f$ has Jacobi
expansion
\begin{equation}\label{rieszjacobipot}
({\mathcal L}^{\alpha,\beta})^{- \nu/2}f(x) = \sum _{k=0}^\infty
\frac{\langle f,P_{k}^{\left( \alpha,\beta\right) }\rangle}{
\hat{h_{k}}^{\left( \alpha ,\beta \right)}}\lambda_k^{-\nu/2}
P_{k}^{\left( \alpha,\beta\right) }(x).
\end{equation}
Now, since
\begin{equation*}
\frac{d}{dx}\left\{ P_{k}^{\left( \alpha,\beta\right) }\left( x\right) \right\} =\frac{\left( k+\alpha +\beta +1\right) }{2%
} P_{k-1}^{\left( \alpha +1,\beta +1\right) }\left(x\right),
\end{equation*}
see Szeg{\"o} \cite{sz},(4.21.7), and using (\ref{rieszjacobipot}) we get that the Jacobi-Riesz transform for  $f\in
L^{2}\left(\left[-1,1\right],\mu_{(\alpha,\beta)}\right)$
 has Jacobi expansion

\begin{equation}\label{RieszJacL2Exp}
R^{\alpha, \beta} f(x)=\sum\limits_{k=1}^{\infty }\frac{\langle
f,P_{k}^{\left( \alpha,\beta\right) } \rangle}{
\hat{h_{k}}^{\left( \alpha ,\beta \right)}}\lambda^{-1/2}_{k}
\frac{(k+\alpha+\beta+1)}{2}\sqrt{1-x^{2}}P^{(\alpha+1,\beta+1)}_{k-1}(x),
\end{equation}
where $\lambda_{k}=k(k+\alpha+\beta+1)$.\\

In particular, the {\em Gegenbauer polynomials}
$\{C^{\lambda}_n\}$, $\lambda>-1/2$ are defined as
\begin{eqnarray}\label{GegebDef}
C^{\lambda}_n(x)&=&\frac{\Gamma(\lambda+1/2)\Gamma(n+2\lambda)}{\Gamma(2\lambda)\Gamma(n+\lambda+1/2)}
P^{(\lambda - 1/2,\lambda -1/2)}_{n}(x).
\end{eqnarray}
 Gegenbauer polynomials are orthogonal with respect to the Gegenbauer measure,  $\mu_{\lambda}$ in $(-1,1)$, defined as,
\begin{eqnarray}
\nonumber\mu_{\lambda}(dx)&=&\omega_{\lambda}(x)dx=\chi_{(-1,1)}(x)\frac{\Gamma(2\lambda+1)}{2^{2\lambda}\left[\Gamma(\lambda+1/2)\right]^{2}}(1-x)^{\lambda-1/2}(1+x)^{\lambda-1/2}dx\\
&=&\chi_{(-1,1)}(x)\frac{\lambda2^{2\lambda}\left[\Gamma(\lambda)\right]^{2}}{2\pi\Gamma(2\lambda)}(1-x^{2})^{\lambda-1/2}dx,\\\nonumber
\end{eqnarray}
by Legendre's duplication formula, see \cite{andaskroy} page 22.\\

Taking $\alpha=\beta=\lambda-1/2$ in (\ref{normPn-1}), then from (\ref{GegebDef}) we have,
\begin{eqnarray}\label{normGege}
 \nonumber\|C^{\lambda}_n\|_{2,\lambda}^{2}&=&\frac{\left[\Gamma(\lambda+1/2)\right]^{2}\left[\Gamma(n+2\lambda)\right]^{2}}{\left[\Gamma(2\lambda)\right]^{2}\left[\Gamma(n+\lambda+1/2)\right]^{2}}\frac{\Gamma(2\lambda+1)\left[\Gamma(n+\lambda+1/2)\right]^{2}
}{(2n+2\lambda)\left[ \Gamma(\lambda+1/2)\right]^{2}\Gamma(n+1)\Gamma(n+2\lambda)}\\
\nonumber&=&\frac{\Gamma(2\lambda+1)\Gamma(n+2\lambda)}{\left[\Gamma(2\lambda)\right]^{2}2(n+\lambda)n!}= \frac{\lambda\Gamma(n+2\lambda)}{\Gamma(2\lambda)(n+\lambda)n!}.\\
\end{eqnarray}

Taking $\alpha=\beta=\lambda-1/2,$ in (\ref{proPn-1}), we get

\begin{eqnarray}\label{dervCn-1}
\nonumber&&\frac{d}{dx}\left\{(1-x^{2})^{\lambda+1/2}C_{n-1}^{\lambda+1}\left(
x\right)\right\}\hspace{7cm}\\
\nonumber&&\hspace{1cm}=-2n\frac{\Gamma(\lambda+3/2)\Gamma(n+2\lambda+1)}{\Gamma(2\lambda+2)\Gamma(n+\lambda+1/2)}\frac{\Gamma(2\lambda)\Gamma(n+\lambda+1/2)}{\Gamma(\lambda+1/2)\Gamma(n+2\lambda)}(1-x^{2})^{\lambda-1/2}C_{n}^{\lambda}\left(x\right)\\
\nonumber&&\hspace{1cm}=-n\frac{(2\lambda+1)(n+2\lambda)}{(2\lambda+1)(2\lambda)}(1-x^{2})^{\lambda-1/2}C_{n}^{\lambda}\left(x\right)=-n\frac{(n+2\lambda)}{2\lambda}(1-x^{2})^{\lambda-1/2}C_{n}^{\lambda}\left(x\right).\\
\end{eqnarray}

Now, for $f\in
L^{2}\left(\left[-1,1\right],\mu_{\lambda}\right)$, from (\ref{RieszJacL2Exp}) we get that its
{\em Gegenbauer-Riesz transform}  will have the expansion
\begin{equation}
R^{\lambda} f(x) = \frac{1}{2} \sum\limits_{k=1}^{\infty }\langle
f,C_{k}^{\lambda} \rangle d_{k}^{\left( \lambda
\right)} (\frac{k+2 \lambda}{k})^{1/2} \sqrt{1-x^{2}}C^{\lambda+1}_{k-1}(x),
\end{equation}
where $d_{k}^{\left( \lambda \right)}= \frac{4
\Gamma(2\lambda)(k+\lambda) \Gamma(k+1)}{\Gamma(k+2\lambda+1)}$.

The  $L^p$-continuity of the Riesz-Jacobi transform
$R^{\alpha, \beta}$, was proved by Zh. Li \cite{li1} and L. Caffarelli \cite{cafacal} in the case $d=1$. In the case  $d\geq1$ $R_i^{\alpha, \beta}, \; i = , \cdots, d,$ are defined analogously, using partial differentiation in (\ref{rieszJac}), and their $L^p$-continuity  was proved by  A. Nowak and P. Sjogren, \cite{noSj} Theorem 5.1. 
\begin{theo} \label{contRieszJacob}
Assume that $1< p < \infty$ and $\alpha, \beta \in
[-1/2,\infty)^d$. There exists a constant $c_p$ such that
\begin{equation}
\| R_i^{\alpha, \beta} f \|_{p,(\alpha, \beta)} \leq c_p \| f
\|_{p,(\alpha, \beta)}.
\end{equation}
for all $ i = 1, \cdots, d$.\\
\end{theo}
For the particular case of the Gegenbauer polynomials this result was obtained in the one dimensional case by 
B. Muckenhoupt and E. Stein in their seminal paper of 1965, \cite{must}.\\

\item {\em Hermite polynomials}:  The {\em Hermite polynomials}
$\{H_n\}_{n}$, are defined as the orthogonal  polynomials
associated with the Gaussian measure in $\mathbb R$, $ \gamma (dx)
=  \frac{e^{-x^2}}{\sqrt{\pi}} dx ,$ i.e.
\begin{equation} \label{proHerOrt}
\int^{\infty}_{-\infty} H_n(y) H_m(y) \, \gamma (dy) =   2^n
n!\delta_{n,m},
\end{equation}
$n,m =0,1,2, \cdots$, with the {\em normalization}
\begin{equation}
H_{2n+1}(0) = 0, \quad H_{2n}(0)= (-1)^n \frac{(2n)!}{n!}.
\end{equation}
We have
\begin{eqnarray}
H_n^{'}(x) &=& 2n H_{n-1}(x),\\
 H_n^{''}(x) &-& 2xH_n^{'}(x) + 2n H_n(x) = 0, \label{eqHerm}
\end{eqnarray}
thus $H_n$ is an eigenfunction of the one dimensional {\em
Ornstein-Uhlenbeck operator} or harmonic oscillator  operator,
\begin{equation}
L= \frac{1}{2} \frac{d^2}{dx^2} - x \frac{d}{dx},
\end{equation}
associated with the eigenvalue $\lambda_n=-n$. \\

The {\em Gaussian-Riesz transform} can be defined formally, see \cite{noStemp}, as
\begin{equation}\label{rieszgauss}
R^{\gamma} = \frac{d}{dx} (-L)^{-1/2}.
\end{equation}
Therefore, for $ f\in
 L^{2}\left(\mathbb{R},\gamma \right)$ with Hermite expansion
 $$ f = \sum\limits_{k=1}^{\infty} \frac{\langle f,H_{k}\rangle}{2^{k}k!} H_{k}, $$
  its Gaussian-Riesz transform has Hermite expansion 
\begin{eqnarray}
R^{\gamma} f(x) = \sum\limits_{k=1}^{\infty}\frac{ \langle
f,H_{k}\rangle}{2^{k}k!} \sqrt{2k}H_{k-1}(x).
\end{eqnarray}
The $L^p$ continuity of the of the
Gaussian-Riesz transform was proved by B. Muckenhoupt  \cite{mu2} in 1969, in the case $d=1$. In the case  $d\geq1$ $R_i^{\gamma}, \; i = , \cdots, d,$ are defined analogously, using partial differentiation in (\ref{rieszgauss}), and  their $L^p$ continuity has been proved by very different ways, using analytic and probabilistic tools,
by  P. A. Meyer \cite{me3}, R. Gundy \cite{gun}, S.
P\'erez and F. Soria \cite{pe}, G. Pissier  \cite{pi}, C. Guti\'errez \cite{gu} and  W. Urbina \cite{ur1}.
\begin{theo} Assume that $1< p < \infty$. There exists a constant $c_p$ such that
\begin{equation}
\| R_i^{\gamma} f \|_{p,\gamma} \leq c_p \| f \|_{p,\gamma}.
\end{equation}
for all $ i = 1, \cdots, d$.\\
\end{theo}

\item  {\em Laguerre polynomials}: For $\alpha>-1$, the {\em Laguerre polynomials} $\{L^{\alpha} _k\}$  are defined as the
orthogonal  polynomials  associated with the Gamma measure on
$(0,\infty)$, $ \mu_{\alpha} (dx) = \chi_{(0,\infty)}(x)
\frac{x^{\alpha}e^{-x}}{\Gamma(\alpha+1)} dx ,$ i.e.
\begin{equation} \label{proLagrOrt}
\int^{\infty}_{0} L^{\alpha}_n(y) L^{\alpha}_m(y) \,\mu_{\alpha}
(dy) =\dbinom{n+\alpha}{n} \delta_{n,m} =
\frac{\Gamma(n+\alpha+1)}{\Gamma(\alpha+1) n!} \delta_{n,m}, 
\end{equation}
$n,m =0,1,2, \cdots$ We have
\begin{equation} \label{eq:aPri}
 (L_k^{\alpha}(x))' = -L_{k-1}^{\alpha+1}(x).
\end{equation}
\begin{equation}\label{eq:aecdif}
x(L_k^{\alpha}(x))''
+(\alpha+1-x)(L_k^{\alpha}(x))'+kL_k^{\alpha}(x)=0.
\end{equation}
thus $L_k^{\alpha}$ is an eigenfunction of the (one-dimensional)
{\em Laguerre differential operator}
$${\mathcal L}^{\alpha}= x \frac{d^2}{dx^2} +(\alpha +1- x) \frac{d}{dx},$$
associated with the eigenvalue $\lambda_k=-k$ \\

The {\em Laguerre-Riesz transform} can be defined formally, see \cite{noStemp}, as
\begin{equation}\label{rieszlaguerre}
R^{\alpha} = \sqrt{x}\frac{d}{dx} ({\mathcal L}^{\alpha})^{-1/2}.
\end{equation}
Therefore for $ f\in
 L^{2}\left((0,\infty),\mu_{\alpha} \right)$ with Laguerre expansion
 $$ f = \sum\limits_{k=0}^{\infty} \frac{ \Gamma(\alpha+1) k!}{\Gamma(k+\alpha+1)}  \langle f,L_k^{\alpha}\rangle L_k^{\alpha}$$
  its Laguerre-Riesz transform has Laguerre expansion
\begin{eqnarray}
R^{\alpha} f(x)= - \sum\limits_{k=1}^{\infty }\ \frac{
\Gamma(\alpha+1) k! }{\Gamma(k+\alpha+1)} (\sqrt{k})^{-1} \sqrt{x}
\langle f,L_k^{\alpha}\rangle L_{k-1}^{\alpha+1}(x).
\end{eqnarray}

The $L^p$ continuity of the Laguerre-Riesz transform was proved by B. Muckenhoupt \cite{mu3}, for the case $d=1$.   In the case $d\geq 1$ $R_i^{\alpha}, \; i = , \cdots, d,$ are defined analogously, using partial differentiation in (\ref{rieszlaguerre}), and  their $L^p$ continuity was proved by A. Nowak \cite{no}
using Littlewood-Paley theory.

\begin{theo} Assume that $1< p < \infty$ and $\alpha \in [-1/2,\infty)^d$. There exists a constant $c_p$  such that
\begin{equation}
\| R_i^{\alpha} f \|_{p,\alpha} \leq c_p \| f \|_{p,\alpha}.
\end{equation}
for all $ i = 1, \cdots, d,$.
\end{theo}
\item Finally, the {\em asymptotic relations} between   Jacobi
polynomials and other classical orthogonal polynomials (see \cite{sz}, (5.3.4) and (5.6.3)) are the following
\begin{enumerate}
\item[i)] With Hermite polynomials,
\begin{equation}\label{jacoherm}
\lim_{\lambda \rightarrow \infty} \lambda^{-n/2} C^{\lambda}_n
(x/\sqrt{\lambda})= \frac{H_n(x)}{n!},
\end{equation}
\item[ii)] With Laguerre polynomials,
\begin{equation}\label{jacolag}
\lim_{\beta \rightarrow \infty} P_{n}^{\left(\alpha ,\beta
\right)}(1-2x/\beta) = L^{\alpha}_n(x).\\
\end{equation}
\end{enumerate}
Both relations holds uniformly in every closed interval of ${\mathbb R}$.
\end{itemize}
We want to thank Professor Luis A. Caffarelli to share with us the idea that such transferences between the Jacobi case and the Hermite and Laguerre cases should be possible.
\section{Main Results}
We we want to obtain the $L^p$-continuity for the Gaussian-Riesz
transform and for the $L^p$-continuity of the Laguerre-Riesz transform
from the $L^p$-continuity of the Jacobi-Riesz transform, using a transference method based on the asymptotic relations between Jacobi polynomials and Hermite and Laguerre polynomials. 

We will start considering the case $p=2$; more precisely we want to prove

\begin{theo} \label{L2result}The $L^2(\mu_{\alpha, \beta})$ boundedness for the Jacobi-Riesz transform
\begin{equation}
\| R^{\alpha, \beta} f \|_{2,(\alpha, \beta)} \leq C_2 \| f
\|_{2,(\alpha, \beta)}
\end{equation}
implies

\begin{enumerate}
\item[i)] the $L^2(\gamma)$  boundedness for the Gaussian-Riesz transform
\begin{equation}
\| R^{\gamma} f \|_{2,\gamma} \leq C_2 \| f \|_{2,\gamma}.
\end{equation}
and\\
\item[ii)] the $L^2(\mu_\alpha)$  boundedness for the Laguerre-Riesz transform
 \begin{equation}
\| R^{\alpha} f \|_{2,\alpha} \leq C_2 \| f \|_{2,\alpha}.
\end{equation}
\end{enumerate}

\end{theo}

For the proof of this we will need the following technical
result,
 \\
\begin{prop} \label{normrel}(norm relations)

\begin{enumerate}
\item[i)] Let  $f\in L^{2}(\mathbb{R},\gamma)$ and define
$f_{\lambda}(x)=f(\sqrt{\lambda}x)1_{[-1,1]}(x)$, then $ f_{\lambda}\in
L^{2}([-1,1],\mu_{\lambda})$ and
\begin{equation}
\lim_{\lambda\to \infty}\|f_{\lambda}\|_{2,\lambda}=\|f\|_{2,\gamma}
\end{equation}

\item[ii)] Let $f\in L^{2}(\mathbb{R},\mu_{\alpha})$ and define $
f_{\beta}(x)=f\left(\frac{\beta}{2}(1-x)\right)1_{[-1,1]}(x)$, then
$f_{\beta}\in L^{2}([-1,1],\mu_{(\alpha,\beta)})$ and
\begin{equation}
\lim_{\beta\to \infty}\|f_{\beta}\|_{2,(\alpha,\beta)}=\|f\|_{2,\alpha}
\end{equation}
\end{enumerate}
\end{prop}

\dem \\

i) First let us prove that $f_{\lambda}\in
L^{2}([-1,1],\mu_{\lambda}).$ For $\;x\in
[-1,1],\;$  
\begin{equation}\label{exp1}
(1-x^{2})^{\lambda-1/2}\leq e^{-\lambda
x^{2}}e^{\frac{x^{2}}{2}}\leq e^{-\lambda x^{2}}e^{\frac{1}{2}}.
\end{equation} 
Then,
\begin{eqnarray*}
&&\frac{\Gamma(2\lambda+1)}{2^{2\lambda}\left[\Gamma\left(\lambda+1/2\right)\right]^{2}}\int^{1}_{-1}(1-x^{2})^{\lambda-1/2}|f_{\lambda}(x)|^{2}dx\hspace{8cm}\\
&&\hspace{4cm} \leq  e^{\frac{1}{2}}\frac{\Gamma(2\lambda+1)}{2^{2\lambda}\left[\Gamma\left(\lambda+1/2\right)\right]^{2}}\int^{1}_{-1}|f(\sqrt{\lambda}x)|^{2}e^{-\lambda x^{2}}dx.\\
\end{eqnarray*}
Now, making the change of variable $u=\sqrt{\lambda}x$, we have 
\begin{eqnarray*}
&&e^{\frac{1}{2}}\frac{\Gamma(2\lambda+1)}{2^{2\lambda}\left[\Gamma\left(\lambda+1/2\right)\right]^{2}}\int^{1}_{-1}|f(\sqrt{\lambda}x)|^{2}e^{-\lambda
x^{2}}dx\hspace{8cm}\\
&&\hspace{3cm}=
e^{\frac{1}{2}}\frac{\Gamma(2\lambda+1)}{\sqrt{\lambda}2^{2\lambda}\left[\Gamma\left(\lambda+1/2\right)\right]^{2}}\int^{\infty}_{-\infty}|f(u)|^{2}1_{[-\sqrt{\lambda},\sqrt{\lambda}]}e^{-u^{2}}du\\
&&\hspace{3cm} \leq e^{\frac{1}{2}}\frac{\Gamma(2\lambda+1)}{\sqrt{\lambda}2^{2\lambda}\left[\Gamma\left(\lambda+1/2\right)\right]^{2}}\int^{\infty}_{-\infty}|f(u)|^{2}e^{-u^{2}}du\\
&&\hspace{3cm}=e^{\frac{1}{2}}\frac{\Gamma(2\lambda+1)\sqrt{\pi}}{\sqrt{\lambda}2^{2\lambda}\left[\Gamma\left(\lambda+1/2\right)\right]^{2}}\|f\|^{2}_{2,\gamma}<\infty,\\
\end{eqnarray*}
as $f\in L^{2}(\mathbb{R},\gamma).$ Therefore $
f_{\lambda}\in L^{2}([-1,1],\mu_{\lambda}).$ \\

On the other hand, using Legendre's duplication formula,
\begin{eqnarray*}
\left\|f_{\lambda}\right\|^{2}_{2,\lambda}
&=&\frac{\Gamma(2\lambda+1)}{2^{2\lambda}\left[\Gamma\left(\lambda+1/2\right)\right]^{2}}\int^{1}_{-1}(1-x^{2})^{\lambda-1/2}|f_{\lambda}(x)|^{2}dx\\
&=&\frac{\lambda\left[\Gamma(\lambda)\right]^{2}2^{2\lambda}}{2\pi\Gamma\left(2\lambda\right)}\int^{1}_{-1}(1-x^{2})^{\lambda-1/2}|f_{\lambda}(x)|^{2}dx.
\end{eqnarray*}
Making the change of variable $x=\frac{y}{\sqrt{\lambda}}$, we get
\begin{eqnarray*}
\left\|f_{\lambda}\right\|^{2}_{2,\lambda}&= &\frac{\lambda
\left[\Gamma(\lambda)\right]^{2}2^{2\lambda}}{2\pi\Gamma\left(2\lambda\right)}\int^{\sqrt{\lambda}}_{-\sqrt{\lambda}}|f_{\lambda}(\frac{y}{\sqrt{\lambda}})|^{2}
(1-\frac{y^{2}}
{\lambda})^{\lambda-1/2}\frac{dy}{\sqrt{\lambda}}.
\end{eqnarray*}
Now, observe that
 $$\displaystyle\lim_{\lambda\to
\infty}f_{\lambda}(\frac{y}{\sqrt{\lambda}})=\displaystyle\lim_{\lambda\to
\infty}f(y)1_{[-\sqrt{\lambda},\sqrt{\lambda}]}(y)=f(y),$$
Set, for all $y \in {\mathbb R}$
$$g(y,\lambda)=
1_{[-\sqrt{\lambda},\sqrt{\lambda}]}(y)\left(f_{\lambda}(\frac{y}{\sqrt{\lambda}})\right)^{2}(1-\frac{y^{2}}
{\lambda})^{\lambda-1/2},$$ 
then clearly
$$\lim_{\lambda\to \infty}g(y,\lambda)=(f(y))^{2}e^{-y^{2}},$$
thus by (\ref{exp1}) we get
\begin{eqnarray*}
|g(y,\lambda)|&=&1_{[-\sqrt{\lambda},\sqrt{\lambda}]}(y)\left(f_{\lambda}(\frac{y}{\sqrt{\lambda}})\right)^{2}(1-\frac{y^{2}}
{\lambda})^{\lambda-1/2}\\
&\leq
&e^{\frac{1}{2}}\left(f_{\lambda}(\frac{y}{\sqrt{\lambda}})\right)^{2}e^{-y^{2}}=e^{\frac{1}{2}}\left(f(y)\right)^{2}e^{-y^{2}}1_{[-\sqrt{\lambda},\sqrt{\lambda}]}(y)<e^{\frac{1}{2}}\left(f(y)\right)^{2}e^{-y^{2}}.
\end{eqnarray*}
Then, for  $f\in L^{2}(\mathbb{R},\gamma),$ by the dominated convergence theorem, we have
\begin{eqnarray*}
\lim_{\lambda\to
\infty}\int^{\infty}_{-\infty}|f_{\lambda}(\frac{y}{\sqrt{\lambda}})|^{2}
1_{[-\sqrt{\lambda},\sqrt{\lambda}]}(y)(1-\frac{y^{2}}
{\lambda})^{\lambda-1/2}dy
&=&\int^{\infty}_{-\infty}\lim_{\lambda\to \infty}g(y,\lambda)dy\\
&=&\int^{\infty}_{-\infty}(f(y))^{2}e^{-y^{2}}dy.
\end{eqnarray*}
On the other hand, using the identity  $$\displaystyle\lim_{z \to
\infty}z^{b-a}\frac{\Gamma(z+a)}{\Gamma(z+b)}=1$$
and Legendre's  duplication
formula, we get,
\begin{eqnarray*}
\displaystyle\lim_{\lambda \to \infty}\frac{\lambda
\left[\Gamma(\lambda)\right]^{2}2^{2\lambda}}{\lambda^{1/2}2\pi\Gamma\left(2\lambda\right)}&=&\displaystyle\lim_{\lambda
\to \infty}\frac{\lambda
\left[\Gamma(\lambda)\right]^{2}2^{2\lambda}}{2\pi(2\pi)^{-1/2}2^{2\lambda-1/2}\Gamma\left(\lambda\right)\Gamma\left(\lambda+\frac{1}{2}\right)\lambda^{1/2}}\\
&=&\displaystyle\lim_{\lambda \to
\infty}\frac{\lambda\Gamma(\lambda)}{\lambda^{1/2}\sqrt{\pi}\Gamma\left(\lambda+\frac{1}{2}\right)}=\displaystyle\lim_{\lambda \to
\infty}\frac{\lambda^{-1/2}\Gamma(\lambda+1)}{\sqrt{\pi}\Gamma\left(\lambda+\frac{1}{2}\right)}=\frac{1}{\sqrt{\pi}}.
\end{eqnarray*}
Then,
\begin{eqnarray*}
\lim\limits_{\lambda\rightarrow \infty }
\left\|f_{\lambda}\right\|^{2}_{2,\lambda}&= &
\lim\limits_{\lambda\rightarrow \infty }\frac{\lambda
\left[\Gamma(\lambda)\right]^{2}2^{2\lambda}}{2\pi\Gamma\left(2\lambda\right)}\int^{\sqrt{\lambda}}_{-\sqrt{\lambda}}|f_{\lambda}(\frac{y}{\sqrt{\lambda}})|^{2}
(1-\frac{y^{2}}
{\lambda})^{\lambda-1/2}\frac{dy}{\sqrt{\lambda}}\\
&=&\frac{1}{\sqrt{\pi}}\int^{\infty}_{-\infty}|f(y)|^2 e^{-y^2} dy
= \left\|f\right\|^{2}_{2,\gamma}.
\end{eqnarray*}

ii) For the Laguerre case we have,
 \begin{eqnarray*}
&&\frac{1}{2^{\alpha+\beta+1}B(\alpha+1,\beta+1)}\int^{1}_{-1}(1-x)^{\alpha}(1+x)^{\beta}|f_{\beta}(x)|^{2}dx\hspace{7cm}\\
&&\hspace{3cm}=\frac{\Gamma(\alpha+\beta+2)}{2^{\alpha+\beta+1}\Gamma(\alpha+1)\Gamma(\beta+1)}\int^{1}_{-1}(1-x)^{\alpha}(1+x)^{\beta}|f_{\beta}(x)|^{2}dx,\\
\end{eqnarray*}
and using the change of variable $x=1-\frac{2y}{\beta}$, 
\begin{eqnarray*}
&&\frac{\Gamma(\alpha+\beta+2)}{2^{\alpha+\beta+1}\Gamma(\alpha+1)\Gamma(\beta+1)}\frac{2}{\beta}\int^{\beta}_{0}\left(\frac{2y}{\beta}\right)^{\alpha}\left(2-\frac{2y}{\beta}\right)^{\beta}|f_{\beta}(1-\frac{2y}{\beta})|^{2}dy\\
&&\hspace{3cm}=\frac{\Gamma(\alpha+\beta+2)}{\Gamma(\alpha+1)\Gamma(\beta+1)\beta^{\alpha+1}}\int^{\beta}_{0}y^{\alpha}\left(1-\frac{y}{\beta}\right)^{\beta}|f(y)|^{2}dy.\\
\end{eqnarray*}
Then
\begin{eqnarray*}
&&\frac{\Gamma(\alpha+\beta+2)}{\Gamma(\alpha+1)\Gamma(\beta+1)\beta^{\alpha+1}}\int^{\beta}_{0}y^{\alpha}\left(1-\frac{y}{\beta}\right)^{\beta}|f(y)|^{2}dy\hspace{7cm}\\
&&\hspace{3cm}\leq \frac{\Gamma(\alpha+\beta+2)}{\Gamma(\alpha+1)\Gamma(\beta+1)\beta^{\alpha+1}}\int^{\infty}_{0}y^{\alpha}e^{-y}|f(y)|^{2}dy\\
&&\hspace{3cm}=\frac{\Gamma(\alpha+\beta+2)}{\Gamma(\beta+1)\beta^{\alpha+1}}
\left\|f\right\|^{2}_{2,\alpha}<\infty,
\end{eqnarray*}
as $f\in L^{2}(\mathbb{R},\mu_{\alpha}),$ and therefore
$f_{\beta}\in L^{2}([-1,1],\mu_{(\alpha,\beta)}).$\\

On the other hand, making the change of variable $x=1-\frac{2y}{\beta}$
\begin{eqnarray*}
\left\|f_{\beta}\right\|^{2}_{2,(\alpha,\beta)}
&=&\frac{\Gamma(\alpha+\beta+2)}{2^{\alpha+\beta+1}\Gamma(\alpha+1)\Gamma(\beta+1)}\frac{2}{\beta}\int^{\beta}_{0}\left(\frac{2y}{\beta}\right)^{\alpha}\left(2-\frac{2y}{\beta}\right)^{\beta}|f_{\beta}(1-\frac{2y}{\beta})|^{2}dy\\
&=&\frac{\Gamma(\beta+(\alpha+2))}{\beta^{\alpha+2}\Gamma(\alpha+1)\Gamma(\beta)}\int^{\infty}_{0}y^{\alpha}\left(1-\frac{y}{\beta}\right)^{\beta}|f_{\beta}(1-\frac{2y}{\beta})|^{2}1_{[0,\beta]}(y)dy.
\end{eqnarray*}
 Now, for  $\beta$ big enough,
\begin{equation}
\frac{\Gamma(\beta+(\alpha+2))}{\Gamma(\beta)}\simeq
\beta^{\alpha+2},
\end{equation}
and therefore
\begin{eqnarray*}
\left\|f_{\beta}\right\|^{2}_{2,(\alpha,\beta)}
&\simeq &\frac{1}{\Gamma(\alpha+1)}\int^{\infty}_{0}y^{\alpha}\left(1-\frac{y}{\beta}\right)^{\beta}|f_{\beta}(1-\frac{2y}{\beta})|^{2}1_{[0,\beta]}(y)dy,\\
\end{eqnarray*}
 for  $\beta$ big enough. Observe that,
 $$\displaystyle\lim_{\beta\to
\infty}f_{\beta}(1-\frac{2y}{\beta})=\displaystyle\lim_{\beta\to
\infty}f(y)1_{[0,\beta]}(y)=f(y),$$
so let
$$s(y,\beta)=y^{\alpha}\left(1-\frac{y}{\beta}\right)^{\beta}|f_{\beta}(1-\frac{2y}{\beta})|^{2}1_{[0,\beta]}(y),$$
then
 $$\lim_{\beta\to
\infty}s(y,\beta)=y^{\alpha}e^{-y}|f(y)|^{2},$$ 
but as
\begin{eqnarray*}
|s(y,\beta)|
&\leq &y^{\alpha}e^{-y}|f_{\beta}(1-\frac{2y}{\beta})|^{2}1_{[0,\beta]}(y) \leq y^{\alpha}e^{-y}|f(y)|^{2},
\end{eqnarray*}
 then,  by the dominated convergent theorem, we get 
\begin{eqnarray*}
\lim_{\beta\to
\infty}\int^{\infty}_{0}y^{\alpha}\left(1-\frac{y}{\beta}\right)^{\beta}|f_{\beta}(1-\frac{2y}{\beta})|^{2}1_{[0,\beta]}(y)dy
&=&\lim_{\beta\to
\infty}\int^{\infty}_{0}s(y,\beta)dy\\
&=&\int^{\infty}_{0}y^{\alpha}e^{-y}|f(y)|^{2}dy.
\end{eqnarray*}
  Then
\begin{eqnarray*}
\lim\limits_{\beta\rightarrow \infty
}\left\|f_{\beta}\right\|^{2}_{2,(\alpha,\beta)}
&=&\frac{1}{\Gamma(\alpha+1)}\int^{\infty}_{0}y^{\alpha}e^{-y}|f(y)|^{2}dy=\left\|f\right\|^{2}_{2,\alpha}.\, \ep 
\end{eqnarray*}
\\
 The next result is analogous to the previous one but refer to the inner product and it will be crucial in the proof of Theorem \ref{L2result}.

\begin{prop} (inner product relations) With the same notation as in Lemma \ref{normrel},
\begin{enumerate}
\item[i)] Let  $f\in L^{2}(\mathbb{R},\gamma),$ then
\begin{equation}
\displaystyle\lim_{\lambda\to\infty}\langle
f_{\lambda},\lambda^{-k/2}C^{\lambda}_{k}\rangle=\langle
f,\frac{H_{k}}{k!}\rangle
\end{equation}
\item[ii)] Let  $f\in
L^{2}(\mathbb{R},\mu_{\alpha}),$ then
\begin{equation}
\displaystyle\lim_{\beta\to\infty}\langle f_{\beta},P_{k}^{(\alpha,\beta)} \rangle =\langle f,L_{k}^{\alpha}\rangle
\end{equation}

\end{enumerate}
\end{prop}
\dem \\
\begin{enumerate}
\item[i)]  Let us prove first that
\begin{equation}
\displaystyle\lim_{\lambda\to\infty}\left\|f_{\lambda}+\lambda^{-k/2}C^{\lambda}_{k}\right\|^{2}_{2,\lambda}=\left\|f+\frac{H_{k}}{k!}\right\|^{2}_{2,\gamma}.
\end{equation}
In fact, given that,
$$\frac{\lambda^{1/2}\left[\Gamma(\lambda)\right]^{2}2^{2\lambda}}{2\pi\Gamma\left(2\lambda\right)}(1-\frac{y^{2}}
{\lambda})^{\lambda-1/2}\left(f_{\lambda}(\frac{y}{\sqrt{\lambda}})+\lambda^{-k/2}C^{\lambda}_{k}(\frac{y}{\sqrt{\lambda}})\right)^{2}1_{[-\sqrt{\lambda},\sqrt{\lambda}]}\geq
0,$$
 by Fatou's lemma we get  

\begin{eqnarray*}
&&\frac{1}{\sqrt{\pi}}\int^{\infty}_{-\infty}\left(f(y)+\frac{H_{k}}{k!}(y)\right)^{2}
e^{-y^2} dy
\hspace{6cm}\\&&\hspace{1cm}\leq\lim\limits_{\lambda\rightarrow
\infty }
\frac{\lambda^{1/2}\left[\Gamma(\lambda)\right]^{2}2^{2\lambda}}{2\pi\Gamma\left(2\lambda\right)}\int^{\sqrt{\lambda}}_{-\sqrt{\lambda}}
\left(f_{\lambda}(\frac{y}{\sqrt{\lambda}})+\lambda^{-k/2}C^{\lambda}_{k}(\frac{y}{\sqrt{\lambda}})\right)^{2}(1-\frac{y^{2}}
{\lambda})^{\lambda-1/2}dy,
\end{eqnarray*}
thus
\begin{eqnarray*}
\left\|f+\frac{H_{k}}{k!}\right\|^{2}_{2,\gamma}\leq\displaystyle\lim_{\lambda\to\infty}\left\|f_{\lambda}+\lambda^{-k/2}C^{\lambda}_{k}\right\|^{2}_{2,\lambda}.
\end{eqnarray*}
On the other hand, let 
$$G(y,\lambda)=\left(\frac{\lambda^{1/2}\left[\Gamma(\lambda)\right]^{2}2^{2\lambda}}{2\pi\Gamma\left(2\lambda\right)}\right)^{1/2}\left(f_{\lambda}(\frac{y}{\sqrt{\lambda}})+\lambda^{-k/2}C^{\lambda}_{k}(\frac{y}{\sqrt{\lambda}})\right)(1-\frac{y^{2}}
{\lambda})^{-1/4}1_{[-\sqrt{\lambda},\sqrt{\lambda}]}(y).$$
Clearly $G\in L^{2}(\mathbb{R},\gamma),$  and moreover
$$\displaystyle\lim_{\lambda\to\infty}G(y,\lambda)=\frac{1}{\pi^{1/4}}\left(f(y)+\frac{H_{k}}{k!}(y)\right)^{2}.$$
Then,
\begin{eqnarray*}
&&\left\|f_{\lambda}+\lambda^{-k/2}C^{\lambda}_{k}\right\|^{2}_{2,\lambda}\hspace{10cm}\\
&&=\int^{\infty}_{-\infty}
\frac{\lambda^{1/2}\left[\Gamma(\lambda)\right]^{2}2^{2\lambda}}{2\pi\Gamma\left(2\lambda\right)}\left(f_{\lambda}(\frac{y}{\sqrt{\lambda}})+\lambda^{-k/2}C^{\lambda}_{k}(\frac{y}{\sqrt{\lambda}})\right)^{2}(1-\frac{y^{2}}
{\lambda})^{\lambda-1/2}1_{[-\sqrt{\lambda},\sqrt{\lambda}]}dy\\
&&=\int^{\infty}_{-\infty}\left(G(y,\lambda)\right)^{2}(1-\frac{y^{2}}
{\lambda})^{\lambda}dy\leq\int^{\infty}_{-\infty}\left(G(y,\lambda)\right)^{2}e^{-y^{2}}dy=\sqrt{\pi}\left\|G(\cdot,\lambda)\right\|^{2}_{2,\gamma}.
\end{eqnarray*}
Then by the continuity of the $L^2$ norm
\begin{eqnarray*}
\displaystyle\lim_{\lambda\to\infty}\left\|f_{\lambda}+\lambda^{-k/2}C^{\lambda}_{k}\right\|^{2}_{2,\lambda}&\leq &\displaystyle\lim_{\lambda\to\infty}\sqrt{\pi}\left\|G(\cdot,\lambda)\right\|^{2}_{2,\gamma}=\sqrt{\pi}\left\|\displaystyle\lim_{\lambda\to\infty}G(\cdot,\lambda)\right\|^{2}_{2,\gamma}\\
&=&\sqrt{\pi}\left\|\frac{1}{\pi^{1/4}}\left(f+\frac{H_{k}}{k!}\right)\right\|^{2}_{2,\gamma}=\left\|f+\frac{H_{k}}{k!}\right\|^{2}_{2,\gamma}.\\
\end{eqnarray*}
Therefore
$$\displaystyle\lim_{\lambda\to\infty}\left\|f_{\lambda}+\lambda^{-k/2}C^{\lambda}_{k}\right\|^{2}_{2,\lambda}=\left\|f+\frac{H_{k}}{k!}\right\|^{2}_{2,\gamma},$$
and analogously,
$$\displaystyle\lim_{\lambda\to\infty}\left\|f_{\lambda}-\lambda^{-k/2}C^{\lambda}_{k}\right\|^{2}_{2,\lambda}=\left\|f-\frac{H_{k}}{k!}\right\|^{2}_{2,\gamma}.$$
Then using the polarization formula,
\begin{eqnarray*}
 \displaystyle\lim_{\lambda\to\infty}\langle f_{\lambda},\lambda^{-k/2}C^{\lambda}_{k}\rangle
 &=&\displaystyle\lim_{\lambda\to\infty}\frac{1}{4}\left[\left\|f_{\lambda}+\lambda^{-k/2}C^{\lambda}_{k}\right\|^{2}_{2,\lambda}-\left\|f_{\lambda}-\lambda^{-k/2}C^{\lambda}_{k}\right\|^{2}_{2,\lambda}\right]\\
&=&\frac{1}{4}\left[\left\|f+\frac{H_{k}}{k!}\right\|^{2}_{2,\gamma}-\left\|f-\frac{H_{k}}{k!}\right\|^{2}_{2,\gamma}\right]=\langle f,\frac{H_{k}}{k!}\rangle.
\end{eqnarray*}

\item[ii)] Analogously as in part i), let us prove first that
\begin{equation}
\lim_{\beta\to\infty}\left\|f_{\beta}+P_{k}^{(\alpha,\beta)}\right\|^{2}_{2,(\alpha,\beta)}=\left\|f+L_{k}^{\alpha}\right\|^{2}_{2,\alpha}.
\end{equation}

By Fatou's lemma, we have
\begin{eqnarray*}
&&\frac{1}{\Gamma(\alpha+1)}\int^{\infty}_{0}\left(f(y)+L_{k}^{\alpha}(y)\right)^{2}
y^{\alpha}e^{-y} dy,
\hspace{7cm}\\&&\leq\lim\limits_{\beta\to\infty }
\frac{\Gamma(\beta+(\alpha+2))}{\beta^{\alpha+2}\Gamma(\alpha+1)\Gamma(\beta)}\int^{\infty}_{0}\left|f_{\beta}(1-\frac{2y}{\beta})+
P_{k}^{(\alpha,\beta)}(1-\frac{2y}{\beta})\right|^{2}1_{[0,\beta]}(y)y^{\alpha}\left(1-\frac{y}{\beta}\right)^{\beta}dy,
\end{eqnarray*}
Then,
\begin{eqnarray*}
\left\|f+L_{k}^{\alpha}\right\|^{2}_{2,\alpha}\leq\displaystyle\lim_{\beta\to\infty}\left\|f_{\beta}+
P_{k}^{(\alpha,\beta)}\right\|^{2}_{2,(\alpha,\beta)}.
\end{eqnarray*}
On the other hand, set

$$S(y,\beta)=\left(\frac{\Gamma(\beta+(\alpha+2))}{\beta^{\alpha+2}\Gamma(\alpha+1)\Gamma(\beta)}\right)^{1/2}\left(f_{\beta}(1-\frac{2y}{\beta})+
P_{k}^{(\alpha,\beta)}(1-\frac{2y}{\beta})\right)1_{[0,\beta]}(y).$$
It is clear that $S\in L^{2}(\mathbb{R},\mu_{\alpha}),$  and moreover$$\displaystyle\lim_{\beta\to\infty}S(y,\beta)=\left(\frac{1}{\Gamma(\alpha+1)}\right)^{1/2}\left(f(y)+L_{k}^{\alpha}(y)\right).$$
Now,
\begin{eqnarray*}
&&\left\|f_{\beta}+
P_{k}^{(\alpha,\beta)}\right\|^{2}_{2,(\alpha,\beta)}\hspace{10cm}\\
&&=\int^{\infty}_{0}\frac{\Gamma(\beta+(\alpha+2))}{\beta^{\alpha+2}\Gamma(\alpha+1)\Gamma(\beta)}\left|f_{\beta}(1-\frac{2y}{\beta})+
P_{k}^{(\alpha,\beta)}(1-\frac{2y}{\beta})\right|^{2}1_{[0,\beta]}y^{\alpha}\left(1-\frac{y}{\beta}\right)^{\beta}dy\\
&&=\int^{\infty}_{0}\left(S(y,\beta)\right)^{2}y^{\alpha}\left(1-\frac{y}{\beta}\right)^{\beta}dy\leq\int^{\infty}_{-\infty}\left(S(y,\beta)\right)^{2}y^{\alpha}e^{-y}dy=\Gamma(\alpha+1)\left\|S(\cdot,\beta)\right\|^{2}_{2,\alpha},
\end{eqnarray*}
and by the continuity of the $L^2$-norm
$$\displaystyle\lim_{\beta\to\infty}\Gamma(\alpha+1)\left\|S(\cdot,\beta)\right\|^{2}_{2,\alpha}=\Gamma(\alpha+1)\left\|\lim_{\beta\to\infty}S(\cdot,\beta)\right\|^{2}_{2,\alpha}.$$
Thus,
\begin{eqnarray*}
\displaystyle\lim_{\beta\to\infty}\left\|f_{\beta}+
P_{k}^{(\alpha,\beta)}\right\|^{2}_{2,(\alpha,\beta)}&\leq &\displaystyle\lim_{\beta\to\infty}\Gamma(\alpha+1)\left\|S(\cdot,\beta)\right\|^{2}_{2,\alpha}=\Gamma(\alpha+1)\left\|\lim_{\beta\to\infty}S(\cdot,\beta)\right\|^{2}_{2,\alpha}\\
&=&\Gamma(\alpha+1)\left\|S(\cdot,\beta)\right\|^{2}_{2,\alpha}=\left\|f+L_{k}^{\alpha}\right\|^{2}_{2,\alpha}.\\
\end{eqnarray*}
Therefore,
$$\displaystyle\lim_{\beta\to\infty}\left\|f_{\beta}+
P_{k}^{(\alpha,\beta)}\right\|^{2}_{2,(\alpha,\beta)}=\left\|f+L_{k}^{\alpha}\right\|^{2}_{2,\alpha},$$
and analogously,
$$\displaystyle\lim_{\beta\to\infty}\left\|f_{\beta}-
P_{k}^{(\alpha,\beta)}\right\|^{2}_{2,(\alpha,\beta)}=\left\|f-L_{k}^{\alpha}\right\|^{2}_{2,\alpha}.$$
Then by the polarization formula,

\begin{eqnarray*}
 \displaystyle\lim_{\beta\to\infty}\langle f_{\beta},P_{k}^{(\alpha,\beta)} \rangle &=&\displaystyle\lim_{\beta\to\infty}\frac{1}{4}\left[\left\|f_{\beta}+P_{k}^{(\alpha,\beta)}\right\|^{2}_{2,(\alpha,\beta)}-\left\|f_{\beta}-P_{k}^{(\alpha,\beta)}\right\|^{2}_{2,(\alpha,\beta)}\right]\\
&=&\frac{1}{4}\left[\left\|f+L_{k}^{\alpha}\right\|^{2}_{2,\alpha}-\left\|f-L_{k}^{\alpha}\right\|^{2}_{2,\alpha}\right]=\langle f,L_{k}^{\alpha}\rangle. \quad   \ep\\
\end{eqnarray*}
\end{enumerate}
We are ready to prove Theorem \ref{L2result},\\

\dem \\

First of all let us observe that by Parseval's identity, for $f
\in L^{2}\left(\left[-1,1\right],\mu_{\lambda}\right)$
\begin{eqnarray*}
\left\|R^{(\alpha,
\beta)}f\right\|^{2}_{2,(\alpha,\beta)}&=&\sum\limits_{k=1}^{\infty
}|\frac{\langle f,P_{k}^{\alpha,\beta}
\rangle}{\hat{h_{n}}^{\left( \alpha ,\beta
\right)}}|^{2}\frac{(k+\alpha+\beta+1)^{2}}{4\lambda_{k}}\left\|\sqrt{1-x^{2}}P^{(\alpha+1,\beta+1)}_{k-1}\right\|^{2}_{2,(\alpha,\beta)}.
\end{eqnarray*}
Since,
\begin{eqnarray*}
&&\left\|\sqrt{1-x^{2}}P^{(\alpha+1,\beta+1)}_{k-1}\right\|^{2}_{2,(\alpha,\beta)}\hspace{15cm}\\
&&=\frac{4(\alpha+1)(\beta+1)}{2^{\alpha+\beta+3}(\alpha+\beta+3)(\alpha+\beta+2)B(\alpha+2,\beta+2)}\\
&&\quad \quad \quad \quad \quad \quad \quad \quad \quad \quad \times \int^{1}_{-1}(1-x)^{\alpha+1}(1+x)^{\beta+1}\left[P^{(\alpha+1,\beta+1)}_{k-1}(x)\right]^{2}dx \hspace{7.1cm}\\
&&=\frac{4(\alpha+1)(\beta+1)}{(\alpha+\beta+3)(\alpha+\beta+2)}\left\|P^{(\alpha+1,\beta+1)}_{k-1}\right\|^{2}_{2,(\alpha+1,\beta+1)}=\frac{4k }{\left( k+\alpha+\beta+1\right)}\left\|P^{(\alpha,\beta)}_{k}\right\|^{2}_{2,(\alpha,\beta)},\\
\end{eqnarray*}
we get,
\begin{eqnarray*}
\left\|R^{\alpha,
\beta}f\right\|^{2}_{2,(\alpha,\beta)}&=&\sum\limits_{k=1}^{\infty }|\frac{\langle
f,P_{k}^{(\alpha,\beta)} \rangle}{\hat{h_{n}}^{\left( \alpha ,\beta
\right)}}|^{2}\frac{(k+\alpha+\beta+1)^{2}}{4k(k+\alpha+\beta+1)}\frac{4k }{\left( k+\alpha+\beta+1\right)}\left\|P^{(\alpha,\beta)}_{k}\right\|^{2}_{2,(\alpha,\beta)}\\
&=&\sum\limits_{k=1}^{\infty }|\langle
f,P_{k}^{(\alpha,\beta)} \rangle|^{2}\\
&&\hspace{0.5cm}\times
\frac{(2k+\alpha+\beta+1)\Gamma(\alpha+1)\Gamma(\beta+1)\Gamma(k+1)\Gamma(k+\alpha+\beta+1)}{\Gamma(\alpha+\beta+2)\Gamma(k+\alpha+1)\Gamma(k+\beta+1)}
\end{eqnarray*}
\begin{enumerate}
\item[i)]  Now, for the Gegenbauer case,
$\alpha=\beta=\lambda-1/2,$  we have
\begin{eqnarray*}
\left\|R^{\lambda}f\right\|^{2}_{2,\lambda}
&=&\sum\limits_{k=1}^{\infty }|\langle
f,C^{\lambda}_{k}\rangle|^{2}
\frac{(2k+2\lambda)\left[\Gamma(\lambda+1/2)\right]^{2}\Gamma(k+1)\Gamma(k+2\lambda)}{\Gamma(2\lambda+1)\left[\Gamma(k+\lambda+1/2)\right]^{2}}\\
&&\quad \quad\quad\quad \quad\quad\quad \quad\quad\quad \quad\quad \quad \quad\quad\times
\frac{\left[\Gamma(2\lambda)\right]^{2}\left[\Gamma(k+\lambda+1/2)\right]^{2}}{\left[\Gamma(\lambda+1/2)\right]^{2}\left[\Gamma(k+2\lambda)\right]^{2}}\\
&=&\sum\limits_{k=1}^{\infty }|\langle
f,C^{\lambda}_{k}\rangle|^{2}
\frac{2(k+\lambda)\left[\Gamma(2\lambda)\right]^{2}\Gamma(k+1)}{\Gamma(2\lambda+1)\Gamma(k+\lambda+1/2)}\\
&=&\sum\limits_{k=1}^{\infty }|\langle
f,C^{\lambda}_{k}\rangle|^{2}
\frac{(k+\lambda)\Gamma(2\lambda)\Gamma(k+1)}{\lambda\Gamma(k+2\lambda)}=\sum\limits_{k=1}^{\infty }|\langle
f,C^{\lambda}_{k}\rangle|^{2} \frac{(k+\lambda)
k!}{\lambda(2\lambda)_{k}}\\
&\geq &\sum\limits_{k=1}^{\infty }|\langle
f,\lambda^{-k/2}C^{\lambda}_{k}\rangle|^{2}
\frac{(\frac{k}{\lambda}+1)k!}{2\sqrt{\pi}(2+1/\lambda)(2+2/\lambda)\ldots (2+(k-1)/\lambda)}.\\
\end{eqnarray*}
Then, we conclude that
\begin{eqnarray*}\label{desig}
&& \sum\limits_{k=1}^{\infty }|\langle
f,\lambda^{-k/2}C^{\lambda}_{k}\rangle|^{2}
\frac{(\frac{k}{\lambda}+1)k!}{2\sqrt{\pi}(2+\frac{1}{\lambda})(2+\frac{2}{\lambda})\ldots
(2+\frac{(k-1)}{\lambda})} \leq \left\|R^{\lambda}
f\right\|^{2}_{2,\lambda}.\hspace{1.5cm}
\end{eqnarray*}

 On the other hand, again using Parseval's identity, we have that the $L^2$-norm
of the Gaussian Riesz transform for  $ f\in
 L^{2}\left(\mathbb{R},\gamma \right)$ is given by
\begin{eqnarray*}
\left\|R^{\gamma}f\right\|^{2}_{2,\gamma}&=&\sum\limits_{k=1}^{\infty }\frac{|\langle
f,H_{k}\rangle|^{2}}{k!2^{k}\sqrt{\pi}}.
\end{eqnarray*}
Now, since  $f_{\lambda}(x)=f(\sqrt{\lambda}x)1_{[-1,1]}(x)\;$  taking  $\lambda \rightarrow \infty$ in
(\ref{desig})
and using the asymptotic relation (\ref{jacoherm}),
we get
\begin{eqnarray*}
\left\|R^{\gamma}f\right\|^{2}_{2,\gamma}&=&\frac{1}{\sqrt{\pi}}\lim\limits_{\lambda\rightarrow
\infty }\sum\limits_{k=1}^{\infty }
\frac{(\frac{k}{\lambda}+1)k!}{2\sqrt{\pi}(2+\frac{1}{\lambda})(2+\frac{2}{\lambda})\ldots
(2+\frac{(k-1)}{\lambda})} |\langle
f_{\lambda},\lambda^{-k/2}C^{\lambda}_{k}\rangle|^{2}\\
&\leq &\frac{1}{\sqrt{\pi}} \lim\limits_{\lambda\rightarrow \infty
}\left\|R^{\lambda}f_{\lambda} \right\|^{2}_{2,\lambda}.
\end{eqnarray*}
Therefore, using Theorem \ref{contRieszJacob} and  Proposition \ref{normrel} i), we get
\begin{eqnarray*}\label{desig4}
 \left\|R^{\gamma}f\right\|^{2}_{2,\gamma}&\leq&
\lim\limits_{\lambda\rightarrow \infty
}\left\|R^{\lambda}f_{\lambda} \right\|^{2}_{2,\lambda} \leq C_2
\lim\limits_{\lambda\rightarrow \infty }\left\|f_{\lambda}
\right\|^{2}_{2,\lambda} =  C_2\left\|f\right\|^2_{2,\gamma}.
\end{eqnarray*}

\item[ii)] Observe that
\begin{equation*}
\| R^{\alpha} f \|_{2,\alpha}^2=  \sum\limits_{k=1}^{\infty }\
\frac{ [\Gamma(\alpha+1) k!]^2 |\langle
f,L_k^{\alpha}\rangle|^2}{k [\Gamma(k+\alpha+1)]^2}  \|\sqrt{x}
L_{k-1}^{\alpha+1}\|_{2,\alpha}^2.
\end{equation*}
But  since
\begin{eqnarray*}
 \|\sqrt{x}  L_{k-1}^{\alpha+1}\|_{2,\alpha}^2 &=& \int_0^{\infty} [L_{k-1}^{\alpha+1}(x)]^2 x^{\alpha} e^{-x} \frac{dx}{\Gamma(\alpha+1)} = \frac{\Gamma(\alpha+2)}{\Gamma(\alpha+1)}  \frac{\Gamma(k+\alpha+1)}{(k-1)! \Gamma(\alpha+2)}\\
 &=&   \frac{\Gamma(k+\alpha+1)}{(k-1)! \Gamma(\alpha+1)},
\end{eqnarray*}
then
\begin{eqnarray*}
\| R^{\alpha} f \|_{2,\alpha}^2&=& \sum\limits_{k=1}^{\infty
}\ \frac{ [\Gamma(\alpha+1) k!]^2 |\langle f,L_k^{\alpha}\rangle|^2}{k [\Gamma(k+\alpha+1)]^2}  \frac{\Gamma(k+\alpha+1)}{(k-1)! \Gamma(\alpha+1)}\\
&=& \sum\limits_{k=1}^{\infty }\ \frac{ \Gamma(\alpha+1) k!
|\langle f,L_k^{\alpha}\rangle|^2}{\Gamma(k+\alpha+1)} .
\end{eqnarray*}
Therefore
$$ \| R^{\alpha} f \|_{2,\alpha}^2 = \lim_{\beta \rightarrow \infty} \sum\limits_{k=1}^{\infty
}\ \frac{ \Gamma(\alpha+1) k! |\langle
f_{\beta},P^{(\alpha,\beta)}_{k}\rangle|^2}{\Gamma(k+\alpha+1)},$$
where, as before
$f_{\beta}(x)=f\left(\frac{\beta}{2}(1-x)\right)1_{[-1,1]}(x).$ Now, for $\beta$ big enough
$$\frac{\Gamma(\beta+(k+\alpha+1))}{\Gamma(\beta+(k+1))\beta^{\alpha}}\simeq 1\hspace{1cm}\mbox{and}\hspace{1cm}\frac{\Gamma(\beta+1)}{\Gamma(\beta+(\alpha+2))\beta^{-\alpha-1}}\simeq 1,$$
then
\begin{eqnarray*}
&&\sum\limits_{k=1}^{\infty }\ |\langle
f_{\beta},P^{(\alpha,\beta)}_{k}\rangle|^2\frac{ \Gamma(\alpha+1)
k!
}{\Gamma(k+\alpha+1)}\hspace{10cm}\\
&&=\sum\limits_{k=1}^{\infty }\ |\langle
f_{\beta},P^{(\alpha,\beta)}_{k}\rangle|^2\frac{ \Gamma(\alpha+1)
k!
}{\Gamma(k+\alpha+1)}\frac{\Gamma(\beta+(k+\alpha+1))}{\Gamma(\beta+(k+1))\beta^{\alpha}}\frac{\Gamma(\beta+1)}{\Gamma(\beta+\alpha+2)\beta^{-\alpha-1}}\\
&&\leq \sum\limits_{k=1}^{\infty }\ |\langle
f_{\beta},P^{(\alpha,\beta)}_{k}\rangle|^2\frac{(2k+\alpha+\beta+1)
k!\Gamma(\alpha+1)\Gamma(\beta+1)\Gamma(k+\beta+\alpha+1)}{\Gamma(\beta+\alpha+2)\Gamma(k+\alpha+1)\Gamma(k+\beta+1)},\\
\end{eqnarray*}
for  $\beta$ big enough. Thus 
\begin{equation}\label{ul}
 \sum\limits_{k=1}^{\infty }\
|\langle f_{\beta},P^{(\alpha,\beta)}_{k}\rangle|^2\frac{
\Gamma(\alpha+1) k! }{\Gamma(k+\alpha+1)}\leq \left\|R^{(\alpha,
\beta)}f_{\beta}\right\|^{2}_{2,(\alpha,\beta)}
\end{equation}
for  $\beta$ big enough. Now, taking $\beta \rightarrow \infty$
and using the asymptotic relation (\ref{jacolag}) in (\ref{ul}) we get,
\begin{eqnarray*}
\| R^{\alpha} f \|_{2,\alpha}^2 &=&\lim_{\beta \rightarrow \infty}
\sum\limits_{k=1}^{\infty }\ \frac{ \Gamma(\alpha+1) k!
|\langle f_{\beta},P^{(\alpha,\beta)}_{k}\rangle|^2}{\Gamma(k+\alpha+1)}\\
&\leq &\lim_{\beta \rightarrow \infty}\left\|R^{(\alpha,
\beta)}f_{\beta}\right\|^{2}_{2,(\alpha,\beta)}.
\end{eqnarray*}
Therefore, using Theorem \ref{contRieszJacob} and  Proposition \ref{normrel} ii), we get
\begin{eqnarray*}\label{desig3}
 \| R^{\alpha} f \|_{2,\alpha}^2 &\leq&
\lim_{\beta \rightarrow \infty}\left\|R^{(\alpha,
\beta)}f_{\beta}\right\|^{2}_{2,(\alpha,\beta)} \leq C_2
\lim\limits_{\lambda\rightarrow \infty }\left\|f_{\beta}
\right\|^{2}_{2,(\alpha,\beta)}=
C_2\left\|f\right\|_{2,\alpha}^2.\ep \\
\end{eqnarray*}
\end{enumerate}

Now we are going to consider the general case $p \neq 2$. For the proof we will follow the argument given by  Betancour et al in \cite{Bet}. \\
\begin{theo}\label{Lpresult}
Let $\alpha,\beta>-1 $ and $1<p<\infty$, then the $L^p(\mu_{\alpha, \beta})$ boundedness for the Jacobi-Riesz transform
\begin{equation}
\| R^{\alpha, \beta} f \|_{p,(\alpha, \beta)} \leq C_p \| f
\|_{p,(\alpha, \beta)}
\end{equation}
implies

\begin{enumerate}
\item[i)] the $L^p(\gamma)$-boundedness for the Gaussian-Riesz transform
\begin{equation}
\| R^{\gamma} f \|_{p,\gamma} \leq C_p \| f \|_{p,\gamma}.
\end{equation}
and\\
\item[ii)] the $L^p(\mu_\alpha)$-boundedness for the Laguerre-Riesz transform
 \begin{equation}
\| R^{\alpha} f \|_{p,\alpha} \leq C_p \| f \|_{p,\alpha}.
\end{equation}
\end{enumerate}

\end{theo}

\dem \\

i) Assume that the operator $R^{\lambda}$ is bounded in
$L^{p}\left([-1,1]\mu_{\lambda}\right)$. Let $\phi\in
C^{\infty}_{0}(\mathbb{R}),$ for each $\lambda>0$ define the function
 $$\phi_{\lambda}(x)=\phi(\sqrt{\lambda}x),\;$$ 
$x\in \mathbb{R}.\;$ For $\lambda$ big enough 
$sop\,\phi_{\lambda}$ is contained in
$[-1,1].$ In what follows  $\;\lambda\;$ will be taken satisfying that condition. \\

Now, from the boundedness of $R^{\lambda}$ we have
\begin{eqnarray*}
\|R^{\lambda}\phi_{\lambda}\|_{L^{p}\left([-1,1]\mu_{\lambda}\right)}\leq
C\|\phi_{\lambda}\|_{L^{p}\left([-1,1]\mu_{\lambda}\right)},
\end{eqnarray*}
in other words,
\begin{eqnarray*}
\left\|\sum\limits_{n=1}^{\infty }\widehat{\phi_{\lambda}}(n)
r_{n}^{\left( \lambda
\right)}\sqrt{1-x^{2}}C^{\lambda+1}_{n-1}\right\|_{L^{p}\left([-1,1]\mu_{\lambda}\right)}\leq
C\|\phi_{\lambda}\|_{L^{p}\left([-1,1]\mu_{\lambda}\right)},
\end{eqnarray*}
where $r_{n}^{\left( \lambda \right)}=d_{n}^{\left( \lambda
\right)}\nu^{-1/2}_{n} \frac{(n+2\lambda)}{2},\;\;$ $\nu_n =
n(n+2\lambda)$ and $d_{n}^{\left( \lambda \right)} = \frac{4
\Gamma(2\lambda)(n+\lambda)\Gamma(n+1)}{\Gamma(n+2\lambda+1)}$.\\

Making the change of variable $x=\frac{y}{\sqrt{\lambda}}$ and taking
$Z(\lambda)=\frac{\lambda^{1/2}\left[\Gamma(\lambda)\right]^{2}2^{2\lambda}}{2\pi\Gamma\left(2\lambda\right)}$
we get
\begin{eqnarray*}
\left\{\int_{-\sqrt{\lambda}}^{\sqrt{\lambda}}\left|\sum\limits_{n=1}^{\infty
}\widehat{\phi_{\lambda}}(n)r_{n}^{\left( \lambda
\right)}\sqrt{1-\frac{y^{2}}{\lambda}}C^{\lambda+1}_{n-1}(\frac{y}{\sqrt{\lambda}})\right|^{p}Z(\lambda)(1-\frac{y^{2}}{\lambda})^{\lambda-1/2}dy\right\}^{1/p}\hspace{3cm}
\end{eqnarray*}
$$ \leq
C\|\phi_{\lambda}\|_{L^{p}\left([-1,1]\mu_{\lambda}\right)},$$ 
and therefore,
\begin{eqnarray*}
\left\{\int_{-\sqrt{\lambda}}^{\sqrt{\lambda}}\left|\sum\limits_{n=1}^{\infty
}\widehat{\phi_{\lambda}}(n)r_{n}^{\left( \lambda
\right)}\left(1-\frac{y^{2}}{\lambda}\right)^{\lambda/p-1/2p+1/2}e^{\frac{y^{2}}{p}}C^{\lambda+1}_{n-1}(\frac{y}{\sqrt{\lambda}})\right|^{p}\frac{e^{-y^{2}}}{\sqrt{\pi}}dy\right\}^{1/p}\hspace{3cm}
\end{eqnarray*}
$$\quad \quad \quad \quad \quad \quad\quad \quad \quad \quad \quad \quad \leq
C(Z(\lambda))^{-1/p}\|\phi_{\lambda}\|_{L^{p}\left([-1,1]\mu_{\lambda}\right)}.$$
On the other hand, we also have
\begin{eqnarray*}
\left\{\int_{-\sqrt{\lambda}}^{\sqrt{\lambda}}\left|\sum\limits_{n=1}^{\infty
}\widehat{\phi_{\lambda}}(n)r_{n}^{\left( \lambda
\right)}\left(1-\frac{y^{2}}{\lambda}\right)^{\lambda/2-1/4+1/2}e^{\frac{y^{2}}{2}}C^{\lambda+1}_{n-1}(\frac{y}{\sqrt{\lambda}})\right|^{2}\frac{e^{-y^{2}}}{\sqrt{\pi}}dy\right\}^{1/2}\hspace{3cm}
\end{eqnarray*}
$$\quad \quad \quad \quad \quad \quad\quad \quad \quad \quad \quad \quad \leq
C(Z(\lambda))^{-1/2}\|\phi_{\lambda}\|_{L^{2}\left([-1,1]\mu_{\lambda}\right)}.$$
Define, for every  $k\in\mathbb{N}\;\;$  and  every $\lambda>0$ such that $\sqrt{\lambda}>k\;$ 
\begin{eqnarray*}
F_{\lambda ,k}(y)=\left\{\begin{array}{lcl}
  \sum\limits_{n=1}^{\infty
}\widehat{\phi_{\lambda}}(n)r_{n}^{\left( \lambda
\right)}C^{\lambda+1}_{n-1}(\frac{y}{\sqrt{\lambda}})\left(1-\frac{y^{2}}{\lambda}\right)^{\lambda/2-1/4+1/2}e^{\frac{y^{2}}{2}}& \mbox{if} & |y|\leq k\\
   \\
  0&\mbox{if} & |y|>k,\\
\end{array}
\right.
\end{eqnarray*}
and
\begin{eqnarray*}
f_{\lambda ,k}(y)=\left\{\begin{array}{lcl}
  \sum\limits_{n=1}^{\infty
}\widehat{\phi_{\lambda}}(n)r_{n}^{\left( \lambda
\right)}C^{\lambda+1}_{n-1}(\frac{y}{\sqrt{\lambda}})\left(1-\frac{y^{2}}{\lambda}\right)^{\lambda/p-1/2p+1/2}e^{\frac{y^{2}}{p}}& \mbox{if} & |y|\leq k\\
   \\
  0&\mbox{if} & |y|>k.\\
\end{array}
\right.
\end{eqnarray*}

From the previous inequalities both series converges for all $y\in
[-\sqrt{\lambda},\sqrt{\lambda}]$, and
$F_{\lambda , k}=f_{\lambda, k}\Omega_{\lambda},\;$ where
$$\Omega_{\lambda}(y)=e^{\frac{y^{2}}{2}-\frac{y^{2}}{p}}\left(1-\frac{y^{2}}{\lambda}\right)^{\lambda/2-\lambda/p-1/4+1/2p}$$
for all  $k\in\mathbb{N}\;\;
 $ and $\sqrt{\lambda}>k.\;$ Moreover
\begin{eqnarray*}
|\Omega_{\lambda}(y)|&=&e^{\frac{y^{2}}{2}-\frac{y^{2}}{p}}\left(1-\frac{y^{2}}{\lambda}\right)^{\lambda/2-\lambda/p-1/4+1/2p}\\
&\leq
&e^{\frac{y^{2}}{2}-\frac{y^{2}}{p}}e^{-\frac{y^{2}}{\lambda}(\lambda/2-\lambda/p-1/4+1/2p)}=e^{y^{2}(\frac{1}{4\lambda}-\frac{1}{2p\lambda})}.\\
\end{eqnarray*}
Therefore, if $p\leq 2$ we have, $|\Omega_{\lambda}(y)|\leq 1$ 
and for  $p>2$, $|\Omega_{\lambda}(y)|\leq
e^{k^{2}(\frac{1}{4\lambda}-\frac{1}{2p\lambda})}.$
 Thus 
$\Omega_{\lambda}$ is bounded in $[-k,k]$. \\

Now
\begin{eqnarray*}
(Z(\lambda))^{-1/p}\|\phi_{\lambda}\|_{L^{p}\left([-1,1]\mu_{\lambda}\right)}&=&(Z(\lambda))^{-1/p}\left\{\int_{-1}^{1}|\phi_{\lambda}(x)|^{p}\frac{\lambda\left[\Gamma(\lambda)\right]^{2}2^{2\lambda}}{2\pi\Gamma\left(2\lambda\right)dx}
\left(1-x^{2}\right)^{\lambda-1/2}\right\}^{1/p}
\end{eqnarray*}
and therefore, by the change of variable $x=\frac{y}{\sqrt{\lambda}}$, we get
\begin{eqnarray}\label{desLpfi}
\nonumber (Z(\lambda))^{-1/p}\|\phi_{\lambda}\|_{L^{p}\left([-1,1]\mu_{\lambda}\right)}&=&(Z(\lambda))^{-1/p}\left\{\int_{-\sqrt{\lambda}}^{\sqrt{\lambda}}|\phi_{\lambda}(\frac{y}{\sqrt{\lambda}})|^{p}Z(\lambda)
\left(1-\frac{y^{2}}{\lambda}\right)^{\lambda-1/2}dy\right\}^{1/p}\\
\nonumber &=&\left\{\int_{-\sqrt{\lambda}}^{\sqrt{\lambda}}|\phi(y)|^{p}\left(1-\frac{y^{2}}{\lambda}\right)^{\lambda-1/2}dy\right\}^{1/p}\\
&\leq
&C\left\{\int_{-\sqrt{\lambda}}^{\sqrt{\lambda}}|\phi(y)|^{p}\frac{e^{-y^{2}}}{\sqrt{\pi}}dy\right\}^{1/p} \leq
C\|\phi\|_{L^{p}\left(\mathbb{R},\gamma\right)},
\end{eqnarray}
and then,
\begin{eqnarray*}
\lim_{\lambda\to
\infty}(Z(\lambda))^{-1/p}\|\phi_{\lambda}\|_{L^{p}\left([-1,1]\mu_{\lambda}\right)}&\leq
& \lim_{\lambda\to
\infty}C\left\{\int_{-\sqrt{\lambda}}^{\sqrt{\lambda}}|\phi(y)|^{p}\frac{e^{-y^{2}}}{\sqrt{\pi}}dy\right\}^{1/p}\\
&=&C\left\{\int_{-\infty}^{\infty}|\phi(y)|^{p}\frac{e^{-y^{2}}}{\sqrt{\pi}}dy\right\}^{1/p} =C\|\phi\|_{L^{p}\left(\mathbb{R},\gamma\right)}.
\end{eqnarray*}

 On the other hand,
\begin{eqnarray}\label{norF1gaus}
\nonumber \|F_{\lambda
,k}\|_{L^{2}\left(\mathbb{R},\gamma\right)}
&=&\left\{\int_{-k}^{k}|F_{\lambda
,k}|^{2}\frac{e^{-y^{2}}}{\sqrt{\pi}}dy\right\}^{1/2}\\
&\leq
&\left\{\int_{-\sqrt{\lambda}}^{\sqrt{\lambda}}|F_{\lambda
,k}|^{2}\frac{e^{-y^{2}}}{\sqrt{\pi}}dy\right\}^{1/2} \leq
C(Z(\lambda))^{-1/2}\|\phi_{\lambda}\|_{L^{2}\left([-1,1]\mu_{\lambda}\right)}.
\end{eqnarray}
Thus, from (\ref{desLpfi}) and (\ref{norF1gaus}),
we have
\begin{eqnarray*}
\|F_{\lambda ,k}\|_{L^{2}\left(\mathbb{R},\gamma\right)}\leq
C\|\phi\|_{L^{2}\left(\mathbb{R},\gamma\right)}.
\end{eqnarray*}
Similarly, as
\begin{eqnarray*}
\nonumber \|F_{\lambda,k}\|_{L^{p}\left(\mathbb{R},\gamma\right)}&=&\left\{\int_{\mathbb{R}}|F_{\lambda
,k}|^{p}\frac{e^{-y^{2}}}{\sqrt{\pi}}dy\right\}^{1/p}=\left\{\int_{-k}^{k}|f_{\lambda
,k}\Omega_{\lambda}(y)|^{p}\frac{e^{-y^{2}}}{\sqrt{\pi}}dy\right\}^{1/p},\\
\end{eqnarray*}
and using the boundedness of $\Omega_{\lambda}$ in $[-k,k]$, we get
\begin{eqnarray}\label{norf1gaus}
\nonumber \|F_{\lambda ,k}\|_{L^{p}\left(\mathbb{R},\gamma\right)}
\nonumber &\leq &C\left\{\int_{-k}^{k}|f_{\lambda
,k}|^{p}\frac{e^{-y^{2}}}{\sqrt{\pi}}dy\right\}^{1/p}\\
&\leq &C\left\{\int_{-\infty}^{\infty}|f_{\lambda
,k}|^{p}\frac{e^{-y^{2}}}{\sqrt{\pi}}dy\right\}^{1/p}\leq 
C(Z(\lambda))^{-1/p}\|\phi_{\lambda}\|_{L^{p}\left([-1,1]\mu_{\lambda}\right)}.
\end{eqnarray}
Then, using (\ref{desLpfi}) and (\ref{norf1gaus}) we get
\begin{eqnarray*}
\|F_{\lambda ,k}\|_{L^{p}\left(\mathbb{R},\gamma\right)}\leq
C\|\phi\|_{L^{p}\left(\mathbb{R},\gamma\right)}
\end{eqnarray*}
for all $\sqrt{\lambda}>k.\;$ Therefore $\{F_{\lambda ,k}\}\;$
is a bounded sequence in
$\;L^{2}\left(\mathbb{R},\gamma\right)\;$ and
$\;L^{p}\left(\mathbb{R},\gamma\right),\;$ thus by the
Bourbaki-Alaoglu's theorem, there exists an increasing sequence
$\{\lambda_{j}\}_{j\in\mathbb{N}},$ such that $\lim_{j\to
\infty}\lambda_{j}=\infty$, and functions
$\;F_{k}\in{L^{2}\left(\mathbb{R},\gamma\right)}$ and
$f_{k}\in{L^{p}\left(\mathbb{R},\gamma\right)}$ satisfying that
\begin{itemize}
\item $ F_{\lambda,k}\to  F_{k},$ as $j\to\infty,$ in the weak topology on ${L^{2}\left(\mathbb{R},\gamma\right)}$
\item $F_{\lambda,k}\to  f_{k},$ as $j\to\infty,$ in the weak topology on ${L^{p}\left(\mathbb{R},\gamma\right)}$
\end{itemize}
Moreover,  $sopF_{k}\cup sopf_{k}\subseteq [-k,k],\;$ and
\begin{eqnarray}\label{AnormFk}
\|F_{k}\|_{L^{2}\left(\mathbb{R},\gamma\right)}\leq\lim_{j\to
\infty}\|F_{\lambda_{j},k}\|_{L^{2}\left(\mathbb{R},\gamma\right)}\leq
C\|\phi\|_{L^{2}\left(\mathbb{R},\gamma\right)},
\end{eqnarray}
and
\begin{eqnarray}\label{Anormfk}
\|f_{k}\|_{L^{p}\left(\mathbb{R},\gamma\right)}\leq
C\|\phi\|_{L^{p}\left(\mathbb{R},\gamma\right)}.
\end{eqnarray}
 Observe that, by Cauchy-Schwartz, defining for $k\in\mathbb{N}$,
$$\tau_k(g)=\int_{-\infty}^{\infty}g(x)\chi_{[-k,k]}(x)dx,$$
$\tau_k\in\left(L^{2}\left(\mathbb{R},\gamma\right)\right)^{\ast}$, and therefore,
\begin{eqnarray*}
\int_{-k}^{k}F_{k}(x)dx&=&\int_{-\infty}^{\infty}F_{k}(x)\chi_{[-k,k]}(x)dx = \tau_k(F_k)\\
&=&\lim_{j\to \infty}  \tau_k(F_{\lambda_{j},k})
=\lim_{j\to
\infty}\int_{-\infty}^{\infty}F_{\lambda_{j},k}(x)\chi_{[-k,k]}(x)dx\\
&=&\int_{-\infty}^{\infty}f_{k}(x)\chi_{[-k,k]}(x)dx = \int_{-k}^{k}f_{k}(x)dx;
\end{eqnarray*}
i.e., $F_{k}=f_{k}\;\;\;a.e\;\;(-k,k),\;$ for all $k$ so
 $F_{k}=f_{k}\;\;\;a.e\;\;\mathbb{R}.\;$
Then from (\ref{Anormfk}), we get
\begin{eqnarray}\label{AnormFk-p}
\|F_{k}\|_{L^{p}\left(\mathbb{R},\gamma\right)}\leq
C\|\phi\|_{L^{p}\left(\mathbb{R},\gamma\right)},
\end{eqnarray}
and therefore, from (\ref{AnormFk}) and (\ref{AnormFk-p}), there exists an increasing sequence $\{\lambda_{j}\}_{j\in
\mathbb{N}}\subset (0,\infty)$, with $\lim_{j\to
\infty}\lambda_{j}=\infty,$, and a function $F\in
L^{p}\left(\mathbb{R},\gamma\right)\cap
L^{2}\left(\mathbb{R},\gamma\right),$ such that
\begin{itemize}
\item For each $k\in \mathbb{N},\;$  $ F_{\lambda_{j},k}\to F,$
as $j\to\infty,$ in the weak topology of
${L^{2}\left(\mathbb{R},\gamma\right)}$ and in the weak topology of ${L^{p}\left(\mathbb{R},\gamma\right)}$, 
and
\item
$\|F\|_{L^{p}\left(\mathbb{R},\gamma\right)}\leq
C\|\phi\|_{L^{p}\left(\mathbb{R},\gamma\right)}.$
\end{itemize}

For each $N\in\mathbb{N},\;\;k\in\mathbb{N}\;\;$ 
and $\lambda$ such that $\sqrt{\lambda}>k,$ let us define
\begin{eqnarray*}
F_{\lambda ,k}^{N}(y)=
\chi_{[-k,k]}(y)\sum\limits_{n=1}^{N}\widehat{\phi_{\lambda}}(n)r_{n}^{\left(
\lambda
\right)}C^{\lambda+1}_{n-1}(\frac{y}{\sqrt{\lambda}})\left(1-\frac{y^{2}}{\lambda}\right)^{\lambda/2-1/4+1/2}e^{\frac{y^{2}}{2}}
\end{eqnarray*}
and
\begin{eqnarray*}
H_{\lambda ,k}^{N}(y)=
\chi_{[-k,k]}(y)\sum\limits_{n=N+1}^{\infty}\widehat{\phi_{\lambda}}(n)r_{n}^{\left(
\lambda
\right)}C^{\lambda+1}_{n-1}(\frac{y}{\sqrt{\lambda}})\left(1-\frac{y^{2}}{\lambda}\right)^{\lambda/2-1/4+1/2}e^{\frac{y^{2}}{2}}.
\end{eqnarray*}
Then, $F_{\lambda ,k}=F_{\lambda ,k}^{N}+H_{\lambda ,k}^{N}.$

Now, we want to prove that for $k\in\mathbb{N}$ and $\lambda>0\;$ such that
$\sqrt{\lambda}>k$ ,
\begin{equation}\label{ineqH}
\int_{-\infty}^{\infty}\left|H_{\lambda
,k}^{N}(y)\right|^{2}\frac{e^{-y^{2}}}{\sqrt{\pi}}dy \leq \frac{C}{N},
\end{equation}
uniformly in $\lambda.$ Take $k\in\mathbb{N}$ and $\lambda>0$; making the change of variable $x = \frac{y^{2}}{\lambda}$ and  using Parseval's identity,
\begin{eqnarray*}
&&\int_{-k}^{k}\left|H_{\lambda
,k}^{N}(y)\right|^{2}\frac{e^{-y^{2}}}{\sqrt{\pi}}dy \hspace{6cm}\\
&&\hspace{1.5cm} \leq
\int_{-\sqrt{\lambda}}^{\sqrt{\lambda}}\left|\sum\limits_{n=N+1}^{\infty}\widehat{\phi_{\lambda}}(n)r_{n}^{\left(
\lambda
\right)}C^{\lambda+1}_{n-1}(\frac{y}{\sqrt{\lambda}})\left(1-\frac{y^{2}}{\lambda}\right)^{\lambda/2-1/4+1/2}e^{\frac{y^{2}}{2}}\right|^{2}\frac{e^{-y^{2}}}{\sqrt{\pi}}dy\\
&&\hspace{1.5cm}=\frac{\sqrt{\lambda}}{\sqrt{\pi}}\int_{-1}^{1}\left|\sum\limits_{n=N+1}^{\infty}\widehat{\phi_{\lambda}}(n)r_{n}^{\left(
\lambda
\right)}C^{\lambda+1}_{n-1}(x)\right|^{2}\left(1-x^{2}\right)^{\lambda+1/2}dx\\
&&\hspace{1.5cm}=\frac{\sqrt{\lambda}}{\sqrt{\pi}}\frac{2\pi\Gamma(2\lambda+2)}{(\lambda+1)\left[\Gamma(\lambda+1)\right]^{2}2^{2\lambda+2}}\sum\limits_{n=N+1}^{\infty}\left|\widehat{\phi_{\lambda}}(n)\right|^{2}\left|r_{n}^{\left(
\lambda\right)}\right|^{2}\left\|C^{\lambda+1}_{n-1}\right\|^{2}_{2,(\lambda+1)}.
\end{eqnarray*}
Using  (\ref{dervCn-1}) we have
$$-\frac{2\lambda}{n(n+2\lambda)}\frac{d}{dx}\left\{(1-x^{2})^{\lambda+1/2}C_{n-1}^{\lambda+1}\left(
x\right)\right\}=C_{n}^{\lambda}\left(x\right)(1-x^{2})^{\lambda-1/2}.$$
Now, integrating by parts, we have
 \begin{eqnarray*}
 \widehat{\phi_{\lambda}}(n)&=&\frac{\lambda\left[\Gamma(\lambda)\right]^{2}2^{2\lambda}}{2\pi\Gamma(2\lambda)}\left\{\frac{\sqrt{\lambda}2\lambda}{n(n+2\lambda)}\int_{-1}^{1}\phi'(\sqrt{\lambda}x)C_{n-1}^{\lambda+1}(x)(1-x^{2})^{\lambda+1/2}dx\right\}\\
&=&\frac{(2\lambda+1)}{(\lambda+1)}\frac{\lambda^{3/2}}{n(n+2\lambda)}\\
&&\quad\quad\quad\quad\quad\times\left\{\frac{(\lambda+1)\left[\Gamma(\lambda+1)\right]^{2}2^{2\lambda+2}}{2\pi\Gamma(2\lambda+2)}\int_{-1}^{1}\phi'(\sqrt{\lambda}x)C_{n-1}^{\lambda+1}(x)(1-x^{2})^{\lambda+1/2}dx\right\}\\
&=&\frac{(2\lambda+1)}{(\lambda+1)}\frac{\lambda^{3/2}}{n(n+2\lambda)}
\left\langle\phi'(\sqrt{\lambda}\cdot),C_{n-1}^{\lambda+1}\right\rangle_{(\lambda+1)}.
\end{eqnarray*}
Then,
\begin{eqnarray*}
&&\int_{-k}^{k}\left|H_{\lambda
,k}^{N}(y)\right|^{2}\frac{e^{-y^{2}}}{\sqrt{\pi}}dy \hspace{6cm}\\
&&\hspace{0.3cm}\leq\frac{\sqrt{\lambda}}{\sqrt{\pi}}\frac{2\pi\Gamma(2\lambda+2)}{(\lambda+1)\left[\Gamma(\lambda+1)\right]^{2}2^{2\lambda+2}}\\
&&\quad\quad\times\sum\limits_{n=N+1}^{\infty}\frac{(2\lambda+1)^{2}}{(\lambda+1)^{2}}\frac{\lambda^{3}}{n^{2}(n+2\lambda)^{2}}\left|\left\langle\phi'(\sqrt{\lambda}\cdot),C_{n-1}^{\lambda+1}\right\rangle_{(\lambda+1)}\right|^{2}\left|r_{n}^{\left(
\lambda\right)}\right|^{2}\left\|C^{\lambda+1}_{n-1}\right\|^{2}_{2,(\lambda+1)}\\
&&\hspace{0.3cm}\leq\frac{\sqrt{\lambda}}{\sqrt{\pi}}\frac{2\pi\Gamma(2\lambda+2)}{(\lambda+1)\left[\Gamma(\lambda+1)\right]^{2}2^{2\lambda+2}}\\
&&\quad\quad\times\sum\limits_{n=N+1}^{\infty}\frac{(2\lambda+1)^{2}}{(\lambda+1)^{2}}\frac{\lambda^{3}}{n^{2}(n+2\lambda)^{2}}\left|\frac{\left\langle\phi'(\sqrt{\lambda}\cdot),C_{n-1}^{\lambda+1}\right\rangle_{(\lambda+1)}}{\left\|C^{\lambda+1}_{n-1}\right\|_{2,(\lambda+1)}}\right|^{2}\left|r_{n}^{\left(
\lambda\right)}\right|^{2}\left\|C^{\lambda+1}_{n-1}\right\|^{4}_{2,(\lambda+1)}.
\end{eqnarray*}
Now, using (\ref{normGege}) and that $\;r_{n}^{\left(
\lambda\right)}= \frac{2
\Gamma(2\lambda)(n+\lambda)\Gamma(n+1)}{\Gamma(n+2\lambda+1)}\left(\frac{(n+2\lambda)}{n}\right)^{1/2}$, we get
\begin{eqnarray*}
\left|r_{n}^{\left(
\lambda\right)}\right|^{2}\left\|C^{\lambda+1}_{n-1}\right\|^{4}_{2,(\lambda+1)}&=&
\frac{4\left[\Gamma(2\lambda)\right]^{2}(n+\lambda)^{2}(n!)^{2}}{\left[\Gamma(n+2\lambda+1)\right]^{2}}\frac{(n+2\lambda)}{n}\\
&&\hspace{1cm}\times\frac{(\lambda+1)^{2}\left[\Gamma(n+2\lambda+1)\right]^{2}n^{2}}{(2\lambda+1)^{2}\left[\Gamma(2\lambda+1)\right]^{2}(n+\lambda)^{2}(n!)^{2}}\\
&=&\frac{4\left[\Gamma(2\lambda)\right]^{2}(n+2\lambda)(\lambda+1)^{2}n}{(2\lambda+1)^{2}\left[\Gamma(2\lambda+1)\right]^{2}}=\frac{(\lambda+1)^{2}(n+2\lambda)n}{(2\lambda+1)^{2}\lambda^{2}}.
\end{eqnarray*}
Therefore,
\begin{eqnarray*}
&&\int_{-k}^{k}\left|H_{\lambda
,k}^{N}(y)\right|^{2}\frac{e^{-y^{2}}}{\sqrt{\pi}}dy\hspace{7cm}\\
&&\hspace{0.5cm}\leq\frac{\sqrt{\lambda}}{\sqrt{\pi}}\frac{2\pi\Gamma(2\lambda+2)}{(\lambda+1)\left[\Gamma(\lambda+1)\right]^{2}2^{2\lambda+2}}
\sum\limits_{n=N+1}^{\infty}\frac{\lambda}{n(n+2\lambda)}\left|\frac{\left\langle\phi'(\sqrt{\lambda}\cdot),C_{n-1}^{\lambda+1}\right\rangle_{(\lambda+1)}}{\left\|C^{\lambda+1}_{n-1}\right\|_{2,(\lambda+1)}}\right|^{2}\\
&&\hspace{0.5cm}\leq\frac{\sqrt{\lambda}}{\sqrt{\pi}}\frac{2\pi\Gamma(2\lambda+2)}{(\lambda+1)\left[\Gamma(\lambda+1)\right]^{2}2^{2\lambda+2}}
\frac{1}{N}\sum\limits_{n=1}^{\infty}\left|\frac{\left\langle\phi'(\sqrt{\lambda}\cdot),C_{n-1}^{\lambda+1}\right\rangle_{(\lambda+1)}}{\left\|C^{\lambda+1}_{n-1}\right\|_{2,(\lambda+1)}}\right|^{2},
\end{eqnarray*}
but,
\begin{eqnarray*}
&&\sum\limits_{n=1}^{\infty}\left|\frac{\left\langle\phi'(\sqrt{\lambda}\cdot),C_{n-1}^{\lambda+1}\right\rangle_{(\lambda+1)}}{\left\|C^{\lambda+1}_{n-1}\right\|_{2,(\lambda+1)}}\right|^{2}
\hspace{7cm}\\&&\hspace{0.7cm}=\frac{(\lambda+1)\left[\Gamma(\lambda+1)\right]^{2}2^{2\lambda+2}}{2\pi\Gamma(2\lambda+2)}\int_{-1}^{1}\left|\phi'(\sqrt{\lambda}x)\right|^{2}(1-x^{2})^{\lambda+1/2}dx\\
&&\hspace{0.7cm}=\frac{(\lambda+1)\left[\Gamma(\lambda+1)\right]^{2}2^{2\lambda+2}\lambda^{-1/2}}{2\pi\Gamma(2\lambda+2)}\int_{-\sqrt{\lambda}}^{\sqrt{\lambda}}\left|\phi'(y)\right|^{2}(1-\frac{y^{2}}{\sqrt{\lambda}})^{\lambda+1/2}dy\\
&&\hspace{0.7cm}\leq C\frac{\sqrt{\pi}(\lambda+1)\left[\Gamma(\lambda+1)\right]^{2}2^{2\lambda+2}\lambda^{-1/2}}{2\pi\Gamma(2\lambda+2)}\int_{-\infty}^{\infty}\left|\phi'(y)\right|^{2} \frac{e^{-y^2}}{\sqrt{\pi}}dy\\
\end{eqnarray*}
Then,
\begin{eqnarray*}
&&\int_{-k}^{k}\left|H_{\lambda
,k}^{N}(y)\right|^{2}\frac{e^{-y^{2}}}{\sqrt{\pi}}dy\leq
\frac{C}{N}\left\|\phi'\right\|^{2}_{L^{2}\left(\mathbb{R},\gamma\right)}
\end{eqnarray*}
i.e.,
\begin{eqnarray*}
\int_{-\infty}^{\infty}\left|H_{\lambda
,k}^{N}(y)\right|^{2}\frac{e^{-y^{2}}}{\sqrt{\pi}}dy&=&\int_{-k}^{k}\left|H_{\lambda
,k}^{N}(y)\right|^{2}\frac{e^{-y^{2}}}{\sqrt{\pi}}dy \leq \frac{C}{N}.
\end{eqnarray*}
Thus  $\{H_{\lambda ,k}^{N}\}_{j\in\mathbb{N}}$ is a bounded sequence on
$\;L^{2}\left(\mathbb{R},\gamma\right),\;$ so by
 Bourbaki-Alaoglu's theorem, there exists a sequence
$\;(\lambda_{j})_{j\in\mathbb{N}},\;$ $\lim_{j\to
\infty}\lambda_{j}=\infty$ such that, for all $N\in \mathbb{N}$,
 $\{H_{\lambda_{j} ,k}^{N}\}_{j\in\mathbb{N}}$ converges weakly in $\;L^{2}\left(\mathbb{R},\gamma\right)\;$ to a function
 $H_{k}^{N}\in L^{2}\left(\mathbb{R},\gamma\right).$ Moreover,
 \begin{equation}\label{norH_n}
\int_{-\infty}^{\infty}\left|H_{k}^{N}(y)\right|^{2}\frac{e^{-y^{2}}}{\sqrt{\pi}}dy\leq \frac{C}{N}.\\
\end{equation}
Then, there exists a non decreasing sequence
$(N_{j})_{j\in\mathbb{N}}$ such that,
\begin{equation}\label{limH_k}
H_{k}^{N_{j}}(y)\longrightarrow 0,\hspace{0.2cm}
a.e.\hspace{0.2cm} y\in\mathbb{R}.
\end{equation}
Defining, for each $m\in\mathbb{N}$,
$F_{k}^{N_{m}}=F-H_{k}^{N_{m}};$ then since $
F_{\lambda,k}\to  F_{k},$ as $j\to\infty$ in the weak topology of ${L^{2}\left(\mathbb{R},\gamma\right)}$ and that $F_{\lambda ,k}^{N}=F_{\lambda ,k}-H_{\lambda ,k}^{N},\;$ we have, for all $m\in\mathbb{N},$
\begin{equation}\label{F_k}
F_{\lambda_{j} ,k}^{N_{m}}\longrightarrow
F_{k}^{N_{m}},\hspace{0.6cm} \mbox{weakly in}
\;\;L^{2}\left(\mathbb{R},\gamma\right)
 \hspace{0.3cm}\mbox{as}\hspace{0.3cm}j\to\infty.
\end{equation}
Also, from (\ref{limH_k}) we get that
\begin{equation}\label{limF_k}
F_{k}^{N_{m}}(y)\longrightarrow F(y),\hspace{0.2cm}
\mbox{as}\; m\to\infty,\;\; a.e.\hspace{0.2cm} y\in\mathbb{R}.
\end{equation}
To finish the proof of i) we need to prove that $F$ is the Hermite-Riesz transform of $\phi$, i.e.
\begin{equation}\label{finalind1}
F(y)=\sum_{n=1}^{\infty}\widehat{\phi}^{\gamma}(n)\frac{1}{2^{n}n!}\sqrt{2n}H_{n-1}(y),
\end{equation}
where $$\widehat{\phi}^{\gamma}(n)=\frac{1}{\sqrt{\pi}}\int_{-\infty}^{\infty}\phi(y)H_{n}(y)e^{-y^{2}}dy,$$
is the $n$-th Fourier-Hermite coefficient of $\phi$. Observe that,
\begin{eqnarray}\label{F_lamnda,k}
\nonumber F_{\lambda ,k}^{N}(y)&=&\sum\limits_{n=1}^{N}\widehat{\phi_{\lambda}}(n)\frac{2\Gamma(2\lambda)(n+\lambda)n!(n+2\lambda)^{1/2}}{\Gamma(n+2\lambda+1)n^{1/2}}
C^{\lambda+1}_{n-1}(\frac{y}{\sqrt{\lambda}})\left(1-\frac{y^{2}}{\lambda}\right)^{\lambda/2+1/4}e^{\frac{y^{2}}{2}}\\
\nonumber &=&\sum\limits_{n=1}^{N}\lambda^{-n/2}\widehat{\phi_{\lambda}}(n)\frac{(\frac{n}{\lambda}+1)n!}{n^{1/2}(\frac{n}{\lambda}+2)^{1/2}2(2+\frac{1}{\lambda})(2+\frac{2}{\lambda})\ldots(2+\frac{n-1}{\lambda})}
\\&&\quad\quad\quad\quad\quad\quad \quad\quad\quad\quad\quad\quad\times 2\lambda^{-\frac{(n-1)}{2}}C^{\lambda+1}_{n-1}(\frac{y}{\sqrt{\lambda}})\left(1-\frac{y^{2}}{\lambda}\right)^{\lambda/2+1/4}e^{\frac{y^{2}}{2}}.\\\nonumber 
\end{eqnarray}
From the asymptotic relation between Jacobi and Hermite polynomials (\ref{jacoherm}), it is easy to see that
\begin{eqnarray*}
\lim_{\lambda \to \infty}\lambda^{-\frac{(n-1)}{2}}C^{\lambda+1}_{n-1}\left(\frac{y}{\sqrt{\lambda}}\right)&=&\lim_{\lambda \to \infty} (\lambda+1)^{-\frac{(n-1)}{2}}\left(1+\frac{1}{\lambda}\right)^{\frac{(n-1)}{2}}C^{\lambda+1}_{n-1}\left(\frac{y}{\sqrt{\lambda+1}}\sqrt{1+\frac{1}{\lambda}}\right)\\
&=& \frac{H_{n-1}(y)}{(n-1)!},
\end{eqnarray*}
uniformly in $[-k,k]$. On the other hand, observe
\begin{eqnarray*}
\lim_{\lambda\to
\infty}\lambda^{-n/2}\widehat{\phi_{\lambda}}(n)&=&\frac{1}{\sqrt{\pi}n!}\int_{-\infty}^{\infty}\phi(y)H_{n}(y)e^{-y^{2}}dy=\frac{1}{n!}\widehat{\phi}^{\gamma}(n).
\end{eqnarray*}
Then taking limit as $\lambda \rightarrow \infty$ in (\ref{F_lamnda,k}) we get
\begin{eqnarray*}
\lim_{\lambda \to \infty} F_{\lambda,k}^{N}(y)&=&\sum_{n=1}^{N}\frac{1}{n!}\widehat{\phi}^{\gamma}(n)\frac{n!}{n^{1/2}2^{1/2}2^{n}}2\frac{H_{n-1}(y)}{(n-1)!}\\
&=&\sum_{n=1}^{N}\widehat{\phi}^{\gamma}(n)\frac{1}{2^{n}n!}\sqrt{2n}H_{n-1}(y).
\end{eqnarray*}
Hence,
\begin{equation}
F_{k}^{N_{m}}(y) =\sum_{n=1}^{N_{m}}\widehat{\phi}^{\gamma}(n)\frac{1}{2^{n}n!}\sqrt{2n}H_{n-1}(y),
\end{equation}
and from (\ref{limF_k}) we finally get (\ref{finalind1}),
$$
F(y)=\sum_{n=1}^{\infty}\widehat{\phi}^{\gamma}(n)\frac{1}{2^{n}n!}\sqrt{2n}H_{n-1}(y).\\
$$

ii)\hspace{0.2cm}The proof of ii) is essentially analogous to i), so we will give less details in this case. Assume that the operator 
$R^{\alpha,\beta}$ is bounded in
$L^{p}\left([-1,1]\mu_{\alpha,\beta}\right)$. 
Let $\phi\in C^{\infty}_{0}((0,\infty))$ and $\beta>0,$ then we have
\begin{eqnarray*}
\|R^{\alpha,\beta}\phi_{\beta}\|_{L^{p}\left([-1,1]\mu_{(\alpha,\beta)}\right)}\leq
C\|\phi_{\beta}\|_{L^{p}\left([-1,1]\mu_{(\alpha,\beta)}\right)},
\end{eqnarray*}
where $\phi_{\beta}(x)=\phi\left(\frac{\beta}{2}(1-x)\right),$
for $x\in(0,\infty)$, i.e.

\begin{eqnarray*}
\left\{\eta_{\alpha,\beta}\int_{-1}^{1}\left|\sum\limits_{n=1}^{\infty}\frac{\widehat{\phi_{\beta}}(n)}{\hat{h_{n}}^{\left(
\alpha,\beta\right)}}\frac{(n+\alpha+\beta+1)^{1/2}}{2n^{1/2}}\sqrt{1-x^{2}}P_{n-1}^{\left(\alpha+1,\beta+1\right)
}(x)\right|^{p}(1-x)^{\alpha}(1+x)^{\beta}dx\right\}^{1/p}\hspace{3cm}
\end{eqnarray*}
$$\hspace{7cm}\leq
C\|\phi_{\beta}\|_{L^{p}\left([-1,1]\mu_{(\alpha,\beta)}\right)}.$$

Now, making the change of variable $x=1-\frac{2}{\beta}y$
\begin{eqnarray*}
&&\left\{\eta_{\alpha,\beta}\frac{2^{\alpha+\beta+1}}{\beta^{\alpha+1}}\int_{0}^{\beta}\left|\sum\limits_{n=1}^{\infty}\frac{\widehat{\phi_{\beta}}(n)}{\hat{h_{n}}^{\left(
\alpha,\beta\right)}}\frac{(n+\alpha+\beta+1)^{1/2}}{n^{1/2}}P_{n-1}^{\left(\alpha+1,\beta+1\right)
}\left(1-\frac{2}{\beta}y\right)\left(\frac{y}{\beta}\right)^{1/2}\right.\right.\hspace{3cm}\\
&&
\left.\left.\hspace{2cm}\times\left(1-\frac{y}{\beta}\right)^{1/2+\beta/p}e^{\frac{y}{p}}
\right|^{p}y^{\alpha}e^{-y}dx\right\}^{1/p}\leq
C\|\phi_{\beta}\|_{L^{p}\left([-1,1]\mu_{(\alpha,\beta)}\right)},
\end{eqnarray*}
and therefore
\begin{eqnarray*}
&&\left\{\int_{0}^{\beta}\left|\sum\limits_{n=1}^{\infty}\frac{\widehat{\phi_{\beta}}(n)}{\hat{h_{n}}^{\left(
\alpha,\beta\right)}}\frac{(n+\alpha+\beta+1)^{1/2}}{n^{1/2}}P_{n-1}^{\left(\alpha+1,\beta+1\right)
}\left(1-\frac{2}{\beta}y\right)\left(\frac{y}{\beta}\right)^{1/2}\right.\right.\hspace{3cm}\\
&&
\left.\left.\hspace{2cm}\times\left(1-\frac{y}{\beta}\right)^{1/2+\beta/p}e^{\frac{y}{p}}
\right|^{p}y^{\alpha}e^{-y}dx\right\}^{1/p}\leq
C(Z_{\alpha,\beta})^{-1/p}\|\phi_{\beta}\|_{L^{p}\left([-1,1]\mu_{(\alpha,\beta)}\right)}
\end{eqnarray*}
where
$Z_{\alpha,\beta}=\frac{\Gamma(\alpha+\beta+2)}{\Gamma(\alpha+1)\beta^{\alpha+2}\Gamma(\beta)}$.
Analogously we get,
\begin{eqnarray*}
&&\left\{\int_{0}^{\beta}\left|\sum\limits_{n=1}^{\infty}\frac{\widehat{\phi_{\beta}}(n)}{\hat{h_{n}}^{\left(
\alpha,\beta\right)}}\frac{(n+\alpha+\beta+1)^{1/2}}{n^{1/2}}P_{n-1}^{\left(\alpha+1,\beta+1\right)
}\left(1-\frac{2}{\beta}y\right)\left(\frac{y}{\beta}\right)^{1/2}\right.\right.\hspace{3cm}\\
&&
\left.\left.\hspace{2cm}\times\left(1-\frac{y}{\beta}\right)^{1/2+\beta/2}e^{\frac{y}{2}}
\right|^{2}y^{\alpha}e^{-y}dx\right\}^{1/2}\leq
C(Z_{\alpha,\beta})^{-1/2}\|\phi_{\beta}\|_{L^{2}\left([-1,1]\mu_{(\alpha,\beta)}\right)}.
\end{eqnarray*}

Now, for each $k\in\mathbb{N} $ and $\beta>0$, such that $\beta>k,$ define the functions
\begin{eqnarray*}
F_{\beta,k}(y)&=&-\chi_{(0,k)}(y)\sum\limits_{n=1}^{\infty}\frac{\widehat{\phi_{\beta}}(n)}{\hat{h_{n}}^{\left(
\alpha,\beta\right)}}\frac{(n+\alpha+\beta+1)^{1/2}}{n^{1/2}}P_{n-1}^{\left(\alpha+1,\beta+1\right)
}\left(1-\frac{2}{\beta}y\right)\\
&&\hspace{4cm}\times\left(\frac{y}{\beta}\right)^{1/2}\left(1-\frac{y}{\beta}\right)^{1/2+\beta/2}e^{\frac{y}{2}},
\end{eqnarray*}
and
\begin{eqnarray*}
f_{\beta
,k}(y)&=&-\chi_{(0,k)}(y)\sum\limits_{n=1}^{\infty}\frac{\widehat{\phi_{\beta}}(n)}{\hat{h_{n}}^{\left(
\alpha,\beta\right)}}\frac{(n+\alpha+\beta+1)^{1/2}}{n^{1/2}}P_{n-1}^{\left(\alpha+1,\beta+1\right)
}\left(1-\frac{2}{\beta}y\right)\\
&&\hspace{4cm}\times\left(\frac{y}{\beta}\right)^{1/2}\left(1-\frac{y}{\beta}\right)^{1/2+\beta/p}e^{\frac{y}{p}}.
\end{eqnarray*}
From the previous inequalities  both series converge for all $y\in (0,\beta)$,
and $F_{\beta , k}=f_{\beta,
k}\Omega_{\beta},\;$ where
$$\Omega_{\beta}(y)=e^{\frac{y}{2}-\frac{y}{p}}\left(1-\frac{y}{\beta}\right)^{\beta/2-\beta/p}$$
for all  $k\in\mathbb{N}\;\;$ and $\beta>0$. Now as,
\begin{eqnarray*}
|\Omega_{\beta}(y)|&=&e^{\frac{y}{2}-\frac{y}{p}}\left(1-\frac{y}{\beta}\right)^{\beta/2-\beta/p}
\leq
e^{\frac{y}{2}-\frac{y}{p}}(e^{-y/\beta})^{\beta/2-\beta/p}=1
\end{eqnarray*}
we conclude that $\Omega_{\beta}$ is bounded in $(0,k)$. On the other hand,
\begin{eqnarray*}
(Z_{\alpha,\beta})^{-1/p}\|\phi_{\beta}\|_{L^{p}\left([-1,1]\mu_{(\alpha,\beta)}\right)}&=&
(Z_{\alpha,\beta})^{-1/p}\left\{\eta_{\alpha,\beta}\int_{-1}^{1}\left|\widehat{\phi_{\beta}}(x)\right|^{p}(1-x)^{\alpha}(1+x)^{\beta}dx\right\}^{1/p},
\end{eqnarray*}
and making the change of variable $x=1-\frac{2}{\beta}y$ we get
\begin{eqnarray*}
&&(Z_{\alpha,\beta})^{-1/p}\left\{\eta_{\alpha,\beta}\int_{-1}^{1}\left|\phi_{\beta}(x)\right|^{p}(1-x)^{\alpha}(1+x)^{\beta}dx\right\}^{1/p}\hspace{5cm}
\\&&\hspace{0.5cm}=
(Z_{\alpha,\beta})^{-1/p}\left\{\eta_{\alpha,\beta}\frac{2}{\beta}\int_{0}^{\beta}\left|\phi_{\beta}(1-\frac{2}{\beta}y)\right|^{p}\left(\frac{2y}{\beta}\right)^{\alpha}\left(2-\frac{2y}{\beta}\right)^{\beta}dy\right\}^{1/p}\\
&&\hspace{0.5cm}=\left\{\int_{0}^{\beta}\left|\phi(y)\right|^{p}y^{\alpha}\left(1-\frac{y}{\beta}\right)^{\beta}dy\right\}^{1/p}\leq\left\{\int_{0}^{\infty}\left|\phi(y)\right|^{p}y^{\alpha}e^{-y}dy\right\}^{1/p}=\|\phi\|_{{L^{p}\left((0,\infty),\mu_{\alpha}\right)}}.
\end{eqnarray*}
Then,
\begin{eqnarray*}
\lim_{\beta\to
\infty}C(Z_{\alpha,\beta})^{-1/p}\|\phi_{\beta}\|_{{L^{p}\left([-1,1]\mu_{(\alpha,\beta)}\right)}}&\leq
& \lim_{\beta\to \infty}C\left\{\int_{0}^{\infty}\left|\phi(y)\right|^{p}y^{\alpha}e^{-y}dy\right\}^{1/p}\\
&=&C\|\phi\|_{{L^{p}\left((0,\infty),\mu_{\alpha}\right)}}.
\end{eqnarray*}
Moreover,
\begin{eqnarray}\label{desLpfiLag}
(Z_{\alpha,\beta})^{-1/p}\|\phi_{\beta}\|_{{L^{p}\left([-1,1],\mu_{(\alpha,\beta)}\right)}}\leq
C\|\phi\|_{{L^{p}\left((0,\infty),\mu_{\alpha}\right)}}.
\end{eqnarray}
Now,
\begin{eqnarray}\label{norF1gausLag}
\nonumber \|F_{\beta
,k}\|_{{L^{2}\left((0,\infty),\mu_{\alpha}\right)}}&=&\left\{\frac{1}{\Gamma(\alpha+1)}\int_{0}^{k}|F_{\beta
,k}|^{2}y^{\alpha}e^{-y}dy\right\}^{1/2}\\
\nonumber &\leq
&C(Z_{\alpha,\beta})^{-1/2}\|\phi_{\beta}\|_{L^{2}\left([-1,1]\mu_{(\alpha,\beta)}\right)}
\end{eqnarray}
and therefore, from (\ref{desLpfiLag}) for $p=2$ and
(\ref{norF1gausLag}), we get
\begin{eqnarray*}
\|F_{\beta ,k}\|_{{L^{2}\left((0,\infty),\mu_{
\alpha}\right)}}\leq
C\|\phi\|_{{L^{2}\left((0,\infty),\mu_{\alpha}\right)}}.
\end{eqnarray*}
Analogously, using that $ \Omega_{\beta}$ is bounded in $(0,k)$,
\begin{eqnarray}\label{norf1gausLag}
\nonumber \|F_{\beta
,k}\|_{{L^{p}\left((0,\infty),\mu_{\alpha}\right)}}&=&\left\{\frac{1}{\Gamma(\alpha+1)}\int_{0}^{k}|f_{\beta
,k}\Omega_{\beta}(y)|^{p}y^{\alpha}e^{-y}dy\right\}^{1/p}\\
\nonumber&\leq &C\left\{\frac{1}{\Gamma(\alpha+1)}\int_{0}^{k}|f_{\beta
,k}|^{p}y^{\alpha}e^{-y}dy\right\}^{1/p}\\
&\leq
&C(Z_{\alpha,\beta})^{-1/p}\|\phi_{\beta}\|_{L^{p}\left([-1,1]\mu_{(\alpha,\beta)}\right)}.
\end{eqnarray}
Then, from (\ref{desLpfiLag}) and (\ref{norf1gausLag}),
\begin{eqnarray*}
\|F_{\beta ,k}\|_{{L^{p}\left((0,\infty),\mu_{\alpha}\right)}}
\leq C\|\phi\|_{{L^{p}\left((0,\infty),\mu_{\alpha}\right)}},
\end{eqnarray*}
for all $\beta>k$. Therefore $\{F_{\beta,k}\}$ is a bounded subsequence in
$\;{L^{2}\left((0,\infty),\mu_{\alpha}\right)}\;$ with 
$\;{L^{p}\left((0,\infty),\mu_{\alpha}\right)}.$ Thus by the
 Bourbaki-Alaoglu's theorem, there exists an increasing sequence $\{\beta_{j}\}_{j\in\mathbb{N}}\;$ with $\lim_{j\to
\infty}\beta_{j}=\infty$ and functions
$\;F_{k}\in{{L^{2}\left((0,\infty),\mu_{\alpha}\right)}}$ and
$f_{k}\in{{L^{p}\left((0,\infty),\mu_{\alpha}\right)}}$
satisfying that

\begin{itemize}
\item $ F_{\beta,k}\to  F_{k},$ as $j\to\infty,$ in the weak topology of ${L^{2}\left((0,\infty),\mu_{\alpha}\right)}$
\item $F_{\beta,k}\to  f_{k},$ as $j\to\infty,$ in the weak topology of ${L^{p}\left((0,\infty),\mu_{\alpha}\right)}.$
\end{itemize}
Then, as in i), we can conclude that there exists an increasing  sequence
 $(\beta_{j})_{j\in
\mathbb{N}}\subset (0,\infty)$ such that $\lim_{j\to
\infty}\beta_{j}=\infty,$ and a function $F\in
L^{p}\left((0,\infty),\mu_{\alpha}\right)\cap
L^{2}\left((0,\infty),\mu_{\alpha}\right),$ such that 

\begin{itemize}
\item for each $k\in \mathbb{N},\;$  $ F_{\beta_{j},k}\to F,$
as $j\to\infty,$ in the weak topology on
${L^{2}\left((0,\infty),\mu_{\alpha}\right)}$ and in the weak topology on ${L^{p}\left((0,\infty),\mu_{\alpha}\right)}$.
 \item $\|F\|_{L^{p}\left((0,\infty),\mu_{\alpha}\right)}\leq
C\|\phi\|_{L^{p}\left((0,\infty),\mu_{\alpha}\right)}$.
\end{itemize}

For each  $N\in\mathbb{N},\;\;k\in\mathbb{N}\;\;$
and $\beta$ such that $\beta>k,$ let us define
\begin{eqnarray*}
&&F_{\beta,k}^{N}(y)=-\chi_{(0,k)}(y)\sum\limits_{n=1}^{N}\frac{\widehat{\phi_{\beta}}(n)}{\hat{h_{n}}^{\left(
\alpha,\beta\right)}}\frac{(n+\alpha+\beta+1)^{1/2}}{n^{1/2}}P_{n-1}^{\left(\alpha+1,\beta+1\right)
}\left(1-\frac{2}{\beta}y\right)\\
&&\hspace{4cm}\times\left(\frac{y}{\beta}\right)^{1/2}\left(1-\frac{y}{\beta}\right)^{1/2+\beta/2}e^{\frac{y}{2}}
\end{eqnarray*}
and
\begin{eqnarray*}
&&H_{\beta,k}^{N}(y)=
-\chi_{(0,k)}(y)\sum\limits_{n=N+1}^{\infty}\frac{\widehat{\phi_{\beta}}(n)}{\hat{h_{n}}^{\left(
\alpha,\beta\right)}}\frac{(n+\alpha+\beta+1)^{1/2}}{n^{1/2}}P_{n-1}^{\left(\alpha+1,\beta+1\right)
}\left(1-\frac{2}{\beta}y\right)\\
&&\hspace{4cm}\times\left(\frac{y}{\beta}\right)^{1/2}\left(1-\frac{y}{\beta}\right)^{1/2+\beta/2}e^{\frac{y}{2}}.
\end{eqnarray*}
Then, $F_{\beta,k}=F_{\beta,k}^{N}+H_{\beta,k}^{N}.$

Again, we want to prove that for $k\in\mathbb{N}$ and $\lambda>0\;$ such that
$\sqrt{\lambda}>k$ ,
\begin{equation}\label{ineqH2}
\int_{0}^{\infty}\left|H_{\beta
,k}^{N}(y)\right|^{2} y^{\alpha}e^{-y} dy \leq \frac{C}{N}.
\end{equation}
uniformly in $\beta$. Take $k\in\mathbb{N},$ and $\beta$ such that $\beta>k$, then using Parseval's identity, the change of variable $x= 1- \frac{2}{\beta}y$ and (\ref{normPn-1}),
\begin{eqnarray*}
&&\frac{1}{\Gamma(\alpha+1)}\int_{0}^{\infty}|H_{\beta,k}^{N}(y)|^{2}y^{\alpha}e^{-y}dy \\
&&\hspace{0.5cm}=\frac{\beta^{\alpha}}{\Gamma(\alpha+1)}\int_{0}^{\beta}\left|\sum\limits_{n=N+1}^{\infty}\frac{\widehat{\phi_{\beta}}(n)}{\hat{h_{n}}^{\left(
\alpha,\beta\right)}}\frac{(n+\alpha+\beta+1)^{1/2}}{n^{1/2}}P_{n-1}^{\left(\alpha+1,\beta+1\right)
}\left(1-\frac{2}{\beta}y\right)\right|^{2}\\
&&\hspace{8cm}\times
\left(\frac{y}{\beta}\right)^{\alpha+1}\left(1-\frac{y}{\beta}\right)^{\beta+1}dy\\&&\hspace{0.5cm}=\frac{\beta^{\alpha+1}}{2^{\alpha+\beta+3}\Gamma(\alpha+1)}\int_{-1}^{1}\left|\sum\limits_{n=N+1}^{\infty}\frac{\widehat{\phi_{\beta}}(n)}{\hat{h_{n}}^{\left(
\alpha,\beta\right)}}\frac{(n+\alpha+\beta+1)^{1/2}}{n^{1/2}}P_{n-1}^{\left(\alpha+1,\beta+1\right)
}\left(x\right)\right|^{2}\\
&&\hspace{8cm}\times
\left(1-x\right)^{\alpha+1}\left(1+x\right)^{\beta+1}dx\\
&&\hspace{0.5cm}=\frac{\beta^{\alpha+1}B(\alpha+2,\beta+2)}{\Gamma(\alpha+1)}\sum\limits_{n=N+1}^{\infty}\left|\frac{\widehat{\phi_{\beta}}(n)}{\hat{h_{n}}^{\left(
\alpha,\beta\right)}}\right|^{2}\frac{(n+\alpha+\beta+1)}{n}\left\|P_{n-1}^{\left(\alpha+1,\beta+1\right)}\right\|^{2}_{\mu_{(\alpha+1,\beta+1)}}\\
&&\hspace{0.5cm}=\frac{\beta^{\alpha+1}B(\alpha+2,\beta+2)}{\Gamma(\alpha+1)}\sum\limits_{n=N+1}^{\infty}\left|\frac{\widehat{\phi_{\beta}}(n)}{\hat{h_{n}}^{\left(
\alpha,\beta\right)}}\right|^{2}\frac{(n+\alpha+\beta+1)}{n}\\
&&\hspace{4cm}\times\frac{(\alpha+\beta+3)(\alpha+\beta+2)n
}{(\alpha+1)(\beta+1)
 \left( n+\alpha+\beta+1\right)}\left\|P^{(\alpha,\beta)}_{n}\right\|^{2}_{\mu_{(\alpha,\beta)}}\\&&\hspace{0.5cm}=\frac{\beta^{\alpha+1}B(\alpha+2,\beta+2)}{\Gamma(\alpha+1)}\sum\limits_{n=N+1}^{\infty}\frac{\left|\widehat{\phi_{\beta}}(n)\right|^{2}}{\hat{h_{n}}^{\left(
\alpha,\beta\right)}}\frac{(\alpha+\beta+3)(\alpha+\beta+2)}{(\alpha+1)(\beta+1)}.
\end{eqnarray*}
Now, integrating by parts and using (\ref{proPn-1})
\begin{eqnarray*}
\widehat{\phi_{\beta}}(n)
&=&-\frac{\beta}{2^{\alpha+\beta+3}B(\alpha+1,\beta+1)n}\int^{1}_{-1}\phi'\left(\frac{\beta}{2}(1-x)\right)P_{n-1}^{\left(
\alpha+1,\beta+1\right) }\left(x\right)\\&&\hspace{6.5cm}\times(1-x)^{\alpha+1}(1+x)^{\beta+1}dx\\
&=&-\frac{\beta(\alpha+1)(\beta+1)}{2^{\alpha+\beta+3}(\alpha+\beta+3)(\alpha+\beta+2)B(\alpha+2,\beta+2)n}\int^{1}_{-1}\phi'\left(\frac{\beta}{2}(1-x)\right)\\&&\hspace{4.5cm}\times
P_{n-1}^{\left(
\alpha+1,\beta+1\right) }\left(x\right)(1-x)^{\alpha+1}(1+x)^{\beta+1}dx\\
&=&-\frac{\beta(\alpha+1)(\beta+1)}{n(\alpha+\beta+3)(\alpha+\beta+2)}\left\langle\phi'\left(\frac{\beta}{2}(1-\cdot)\right),P_{n-1}^{\left(
\alpha+1,\beta+1\right) }\right\rangle_{(\alpha+1,\beta+1)}.
\end{eqnarray*}
Then,
\begin{eqnarray*}
&&\frac{1}{\Gamma(\alpha+1)}\int_{0}^{k}|H_{\beta,k}^{N}(y)|^{2}y^{\alpha}e^{-y}dy\hspace{8cm}\\
\end{eqnarray*}
\begin{eqnarray*}
&&=\frac{\beta^{\alpha+3}B(\alpha+2,\beta+2)}{\Gamma(\alpha+1)}\\
&&\hspace{1.5cm}\times\sum\limits_{n=N+1}^{\infty}\frac{\left|\left\langle\phi'\left(\frac{\beta}{2}(1-\cdot)\right),P_{n-1}^{\left(
\alpha+1,\beta+1\right)
}\right\rangle_{(\alpha+1,\beta+1)}\right|^{2}}{\left\|P^{(\alpha+1,\beta+1)}_{n-1}\right\|^{2}_{\mu
_{(\alpha+1,\beta+1)}}}\frac{1}{n(n+\alpha+\beta+1)}\\
&&\leq\frac{\beta^{\alpha+2}B(\alpha+2,\beta+2)}{\Gamma(\alpha+1)N}\sum\limits_{n=N+1}^{\infty}\frac{\left|\left\langle\phi'\left(\frac{\beta}{2}(1-\cdot)\right),P_{n-1}^{\left(
\alpha+1,\beta+1\right)
}\right\rangle_{(\alpha+1,\beta+1)}\right|^{2}}{\left\|P^{(\alpha+1,\beta+1)}_{n-1}\right\|^{2}_{\mu_{(\alpha+1,\beta+1)}}}.
\end{eqnarray*}
But, making the change of variable $y= \frac{\beta}{2}(1-x)$, i.e. $x = 1- \frac{2y}{\beta},$
\begin{eqnarray*}
&&\sum\limits_{n=N+1}^{\infty}\frac{\left|\left\langle\phi'\left(\frac{\beta}{2}(1-\cdot)\right),P_{n-1}^{\left(
\alpha+1,\beta+1\right)
}\right\rangle_{(\alpha+1,\beta+1)}\right|^{2}}{\left\|P^{(\alpha+1,\beta+1)}_{n-1}\right\|^{2}_{\mu_{(\alpha+1,\beta+1)}}}\hspace{4cm}\\
&&=\frac{1}{B(\alpha+2,\beta+2)}\frac{\beta}{\beta^{\alpha+2}}\int_{0}^{\beta}\left|\phi'(y)\right|^{2}y^{\alpha+1}\left(1-\frac{y}{\beta}\right)^{\beta+1}dy\\
&&\leq\frac{1}{B(\alpha+2,\beta+2)}\frac{1}{\beta^{\alpha+2}}\int_{0}^{\infty}\left|\phi'(y)\right|^{2}y^{\alpha+1}e^{-y}dy =\frac{\Gamma(\alpha+1)}{B(\alpha+2,\beta+2)\beta^{\alpha+2}}\left\|\phi'\right\|^{2}_{\mu_{(\alpha+1)}}.
\end{eqnarray*}
Then,
\begin{eqnarray*}
\frac{1}{\Gamma(\alpha+1)}\int_{0}^{k}|H_{\beta,k}^{N}(y)|^{2}y^{\alpha}e^{-y}dy&\leq&\frac{\beta^{\alpha+2}B(\alpha+2,\beta+2)}{\Gamma(\alpha+1)N}\frac{\Gamma(\alpha+1)}{B(\alpha+2,\beta+2)\beta^{\alpha+2}}\left\|\phi'\right\|^{2}_{\mu_{(\alpha+1)}}\\
&=&\frac{1}{N}\left\|\phi'\right\|^{2}_{\mu_{(\alpha+1)}},
\end{eqnarray*}
hence,
\begin{eqnarray*}
\frac{1}{\Gamma(\alpha+1)}\int_{0}^{\infty}|H_{\beta,k}^{N}(y)|^{2}y^{\alpha}e^{-y}dy\leq 
\frac{C}{N}\left\|\phi'\right\|^{2}_{\mu_{(\alpha+1)}}.\\
\end{eqnarray*}

Therefore, again by Bourbaki-Alaoglu's theorem, there exists a sequence
$\{\beta_{j}\}_{j\in\mathbb{N}}$ with $\lim_{j\to
\infty}\beta_{j}=\infty$ such that, for all $N\in \mathbb{N}$,
 $\{H_{\beta_{j} ,k}^{N}\}_{j\in\mathbb{N}}$ converges weakly in
 $\;L^{2}\left((0,\infty),\mu_{\alpha}\right)\;$ to a function
 $H_{k}^{N}\in L^{2}\left((0,\infty),\mu_{\alpha}\right).$ Moreover,
 \begin{equation}\label{norH_nLag}
\int_{-\infty}^{\infty}\left|H_{k}^{N}(y)\right|^{2}y^{\alpha}\frac{e^{-y}}{\Gamma(\alpha+1)}dy\leq \frac{C}{N}.
\end{equation}
Then, there exists a non decreasing sequence
$\{N_{j}\}_{j\in\mathbb{N}}$ such that
\begin{equation}\label{limH_kLag}
H_{k}^{N_{j}}(y)\longrightarrow 0,\hspace{0.2cm}
a.e.\hspace{0.2cm} y\in\mathbb{R}.
\end{equation}
Now, defining, for each $m\in\mathbb{N}$,
$F_{k}^{N_{m}}=F-H_{k}^{N_{m}},$ then since  $
F_{\beta,k}\to F_{k},$ as $j\to\infty,$ in the weak topology of  $L^{2}\left((0,\infty),\mu_{\alpha}\right)$ and that$F_{\beta ,k}^{N}=F_{\beta ,k}-H_{\beta ,k}^{N},\;$ we have that for all $m\in\mathbb{N},$
\begin{equation}\label{F_kLag}
F_{\beta_{j} ,k}^{N_{m}}\longrightarrow
F_{k}^{N_{m}},\hspace{0.3cm} \mbox{weakly in}
\;\;L^{2}\left((0,\infty),\mu_{\alpha}\right)
 \hspace{0.3cm}\mbox{as}\hspace{0.3cm}j\to\infty.
\end{equation}
Additionally, from (\ref{limH_kLag}) we have
\begin{equation}\label{limF_kLag}
F_{k}^{N_{m}}(y)\longrightarrow F(y),\hspace{0.2cm}
\mbox{as}\; m\to\infty,\;\; a.e.\hspace{0.2cm} y\in\mathbb{R}.
\end{equation}
To finish the proof of ii) we need to prove that F is the Laguerre-Riesz transform of $\phi$, i.e.
\begin{equation}\label{finalind2}
F(y)=-\sum\limits_{n=1}^{\infty}\frac{n!\Gamma(\alpha+1)}{\left(
n+\alpha
+1\right)}\widehat{\phi}^{\alpha}(n)\frac{y^{1/2}}{n^{1/2}}L_{n-1}^{\alpha+1}(y).
\end{equation}

From the asymptotic relation between Jacobi and Laguerre polynomials (\ref{jacolag}), it is easy to see that

\begin{eqnarray*}
\lim_{\beta \to \infty} P^{(\alpha+1,\beta+1)}_{n-1}\left(1-\frac{2y}{\beta}\right)&=& \lim_{\beta \to \infty}P^{(\alpha+1,\beta+1)}_{n-1}\left(1-\frac{2y}{\beta+1}\frac{(\beta+1)}{\beta}\right)= L_{n-1}^{\alpha+1}(y).
\end{eqnarray*}

On the other hand,
\begin{eqnarray*}
\lim_{\beta\to
\infty}\widehat{\phi_{\beta}}(n)&=&\frac{1}{\Gamma(\alpha+1)}\int_{0}^{\infty}\phi(y)L_{n}^{\alpha}(y)y^{\alpha}e^{-y}dy= \widehat{\phi}^{\alpha}(n),
\end{eqnarray*}
 and
 \begin{eqnarray*}
\lim_{\beta\to
\infty}\frac{1}{\hat{h_{n}}^{\left(\alpha,\beta\right)}}&=&\lim_{\beta\to\infty}(\frac{2n}{\beta}+\frac{\alpha}{\beta} +1+\frac{1}{\beta})n!\frac{\Gamma(\alpha+1)}{\left( n+\alpha
+1\right)}\frac{\beta^{\alpha+1}\Gamma(\beta+1)}{\Gamma(\alpha+\beta+2)}\frac{\Gamma\left(\beta+n+\alpha+1\right)}{\beta^{\alpha}\Gamma
\left( n+\beta +1\right) }\\
&=&\frac{n!\Gamma(\alpha+1)}{\left( n+\alpha +1\right)}.
\end{eqnarray*}

Then, for $y\in (0,k)$
\begin{eqnarray*}
\lim_{\beta\to \infty}F_{\beta,k}^{N}(y)&=&-\lim_{\beta\to
\infty}\sum\limits_{n=1}^{N}\frac{\widehat{\phi_{\beta}}(n)}{\hat{h_{n}}^{\left(
\alpha,\beta\right)}}\frac{(n+\alpha+\beta+1)^{1/2}}{n^{1/2}}P_{n-1}^{\left(\alpha+1,\beta+1\right)
}\left(1-\frac{2}{\beta}y\right)\\
&&\hspace{4cm}\times\left(\frac{y}{\beta}\right)^{1/2}\left(1-\frac{y}{\beta}\right)^{1/2+\beta/2}e^{\frac{y}{2}}\\
&=&-\lim_{\beta\to
\infty}\sum\limits_{n=1}^{N}\frac{\widehat{\phi_{\beta}}(n)}{\hat{h_{n}}^{\left(
\alpha,\beta\right)}}\frac{(\frac{n}{\beta}+\frac{\alpha}{\beta}+1+\frac{1}{\beta})^{1/2}}{n^{1/2}}P_{n-1}^{\left(\alpha+1,\beta+1\right)
}\left(1-\frac{2}{\beta}y\right)\\
&&\hspace{4cm}\times
y^{1/2}\left(1-\frac{y}{\beta}\right)^{1/2+\beta/2}e^{\frac{y}{2}}\\
&=&-\sum\limits_{n=1}^{N}\frac{n!\Gamma(\alpha+1)}{\left( n+\alpha
+1\right)}\widehat{\phi}^{\alpha}(n)\frac{y^{1/2}}{n^{1/2}}L_{n-1}^{\alpha+1}(y).
\end{eqnarray*}

Thus,
\begin{equation}
F_{k}^{N_{m}}=-\sum\limits_{n=1}^{N_{m}}\frac{n!\Gamma(\alpha+1)}{\left(
n+\alpha
+1\right)}\widehat{\phi}^{\alpha}(n)\frac{y^{1/2}}{n^{1/2}}L_{n-1}^{\alpha+1}(y),
\end{equation}
and therefore we get (\ref{finalind2}),
$$
F(y)=-\sum\limits_{n=1}^{\infty}\frac{n!\Gamma(\alpha+1)}{\left(
n+\alpha
+1\right)}\widehat{\phi}^{\alpha}(n)\frac{y^{1/2}}{n^{1/2}}L_{n-1}^{\alpha+1}(y). \; \ep
$$

\end{document}